\documentclass[10pt]{amsart}
\usepackage{rotating}
\usepackage{url}
\usepackage[all]{xy}
\usepackage{color}
\usepackage{multirow}
\usepackage[table]{xcolor}

\usepackage{amssymb,verbatim,url}
\definecolor{g}{rgb}{.75,.75,.75}
\newcommand{\Z}{{\mathbb Z}}

\newcommand{\Q}{{\mathbb Q}}
\newcommand{\F}{{\mathbb F}}
\newcommand{\C}{{\mathbb C}}

\newcommand{\Gal}{\mbox{Gal}}
\newcommand{\ord}{\mbox{ord}}

\newcommand{\z}{\hspace{5pt}}

\newcommand{\LCM}{\mbox{LCM}}
\newcommand{\cmmt}[1]{}
\newtheorem{Theorem}{Theorem}[section]

\newtheorem{Conjecture}[Theorem]{Conjecture}

\allowdisplaybreaks

\author{David P.\ Roberts}
\address{Division of Science and Mathematics, University of
  Minnesota-Morris, Morris, MN 56267}
\email{roberts@morris.umn.edu}

\title{Newforms with rational coefficients}
\begin{document}
\maketitle 
\begin{abstract}
We consider the set of classical newforms with
rational coefficients and no complex multiplication.   
We study the distribution of quadratic-twist classes of these
forms with respect to weight $k$ and minimal level $N$.
We conjecture that for each weight $k \geq 6$, there
are only finitely many classes.   In large weights,
we make this conjecture effective:  in weights $18 \leq k \leq 24$, 
all classes have $N \leq 30$, in weights 
$26 \leq k \leq 50$, all classes have $N \in \{2,6\}$, and
in weights $k \geq 52$, there are no classes at all.  
We study some of the newforms appearing on
our conjecturally complete list 
in more detail, especially in the cases 
$N=2$, $3$, $4$, $6$, and $8$, 
where formulas can be kept nearly
as simple as those for the classical case $N=1$.  
\end{abstract}

\section{Introduction}
\label{intro}

\subsection{A finiteness conjecture, effective in large weights.}
    Classical newforms, as reviewed in
     next section, are certain power series $g=\sum_{n=1}^\infty a_n q^n \in \C[[q]]$ which play an important role in arithmetic geometry.  This paper is a contribution to 
cataloging newforms for which all the coefficients $a_n$ are rational.  
  We exclude newforms with complex
multiplication, as CM newforms with rational coefficients 
have been comprehensively treated by Sch\"utt \cite{Sch09}.  
We quotient
out by the operation of quadratic twisting, thereby replacing infinitely many 
newforms $g^\chi$ by a single twist-class of newforms $[g]$.  Our focus is
then on the finite sets $Q_k(N)$ of twist-classes of newforms with rational 
coefficients, no complex multiplication,
 and minimal level $N$.
Here the weight $k$ runs over positive even integers, while the minimal level
$N$ runs over positive integers. 

\begin{table}[htb]
{\small
\[
{\renewcommand{\arraycolsep}{.95pt}
\def\arraystretch{.92}
\begin{array}{|r|rrrrrrrrrrrrrrrrrrrrrrrrr|}
\hline
 \text{} & \z 2 & \z 4 & \z 6 & \z 8 & \z 10 & \z 12 & 14 & 16 & 18 & 20 &
   22 & 24 & 26 &  28 & 30 & 32 & 34 & 36 & 38 & 40 & 42 & 44 & 46 & 48 & 50  \\
 \hline
 1 & \text{} & \text{} & \text{} & \text{} & \text{} & 1
   & \text{} & 1 & 1 & 1 & 1 & \text{} & 1 & \text{} &&&&&&&&&&&  \\
 2 & \text{} & \text{} & \text{} & 1 & 1 & \text{} & 2 &
   1 & 1 & 2 & 2 & 1 & 1 & 2 & 2 & 1 & 1 & 2 &  & 1 & 1 &  & & 1 &   \\
 3 & \text{} & \text{} & 1 & 1 & 2 & 1 & 1 & 2 & 1 & 1 &
   {\bf \underline{2}} & 1 & \text{} & \text{} &&&&&&&&&&&\\
 4 & \text{} & \text{} & 1 & \text{} & 1 & 1 & 1 & 1 &
   \text{} & 1 & \text{} & \text{} & \text{} & \text{}
 &&&&&&&&&&&  \\
 5 & \text{} & 1 & 1 & 1 & 1 & 1 & \text{} & \text{} &
   \text{} & \text{} & \text{} & \text{} & \text{} &
   \text{}&&&&&&&&&&& \\
 6 & \text{} & 1 & 1 & 1 & 1 & 3 & 1 & 3 & 3 & 3 & 3 & 3
   & 3 & 3 & 3 & 3 & 3 & 1 & 3 & 1 & 1 & 1 & 1 &  & 1  \\
 7 & \text{} & 1 & 1 & 1 & \text{} & \text{} & \text{} &
   \text{} & \text{} & \text{} & \text{} & \text{} &
   \text{} & \text{} &&&&&&&&&&&\\
 8 & \text{} & 1 & 1 & 2 & 2 & 1 & 1 & {\bf \underline{2}} & \text{} &
   \text{} & \text{} & \text{} & \text{} & \text{} &&&&&&&&&&&\\
 9 & \text{} & \text{} & \text{} & \text{} & \text{} &
   \text{} & \text{} & \text{} & \text{} & \text{} &
   \text{} & \text{} & \text{} & \text{}&&&&&&&&&&& \\
 10 & \text{} & 1 & 3 & 1 & 3 & 3 & 3 & 3 & 1 & 3 & 1 &
   \text{} & 1 & \text{}&&&&&&&&&&& \\
 11 & 1 & \text{} & 1 & \text{} & \text{} & \text{} &
   \text{} & \text{} & \text{} & \text{} & \text{} &
   \text{} & \text{} & \text{} &&&&&&&&&&&\\
 12 & \text{} & 1 & \text{} & 2 & 1 & 2 & 2 & 1 & 2 &
   \text{} & 1 & \text{} & \text{} & \text{} &&&&&&&&&&&\\
 13 & \text{} & 1 & \text{} & {\bf {1}} & \text{} & \text{} &
   \text{} & \text{} & \text{} & \text{} & \text{} &
   \text{} & \text{} & \text{} &&&&&&&&&&&\\
 14 & 1 & 2 & 2 & 2 & 2 & 2 & 2 & 1 & 1 & \text{} &
   \text{} & \text{} & \text{} & \text{} &&&&&&&&&&& \\
 15 & 1 & 2 & 2 & 2 & 2 & 1 & 1 & \text{} & \text{} &
   \text{} & {\bf {1}} & \text{} & \text{} & \text{} &&&&&&&&&&& \\
 16 & \text{} & \text{} & \text{} & \text{} & \text{} &
   \text{} & \text{} & \text{} & \text{} & \text{} &
   \text{} & \text{} & \text{} & \text{} &&&&&&&&&&&\\
 17 & 1 & 1 & {\bf \underline{2}} & {\bf {1}} & \text{} & \text{} & \text{} &
   \text{} & \text{} & \text{} & \text{} & \text{} &
   \text{} & \text{} &&&&&&&&&&& \\
 18 & \text{} & \text{} & 1 & \text{} & 1 & 1 & 1 & 1 &
   \text{} & 1 & \text{} & \text{} & \text{} & \text{} &&&&&&&&&&&
   \\
 19 & 1 & 1 & {\bf \underline{2}} & \text{} & \text{} & \text{} & \text{}
   & \text{} & \text{} & \text{} & \text{} & \text{} &
   \text{} & \text{} &&&&&&&&&&& \\
 20 & 1 & 1 & 1 & 1 & 1 & 1 & \text{} & \text{} &
   \text{} & \text{} & \text{} & \text{} & \text{} &
   \text{} &&&&&&&&&&& \\
 21 & 1 & 2 & \mathbf{{\underline{4}}} & 1 & 1 & \mathbf{{2}} & \text{} & \text{} &
   \text{} & \text{} & \text{} & \text{} & \text{} &
   \text{} &&&&&&&&&&& \\
 22 & \text{} & 3 & 3 & 3 & 3 & \text{} & \text{} &
   \text{} & \text{} & {\bf{1}} & \text{} & \text{} & \text{} &
   \text{} &&&&&&&&&&&\\
 23 & \text{} & 1 & \text{} & \text{} & \text{} &
   \text{} & \text{} & \text{} & \text{} & \text{} &
   \text{} & \text{} & \text{} & \text{} &&&&&&&&&&& \\
 24 & 1 & 1 & 3 & 3 & 3 & 3 & 1 & 1 & \text{} & \text{}
   & \text{} & \text{} & \text{} & \text{} &&&&&&&&&&&\\
 25 & \text{} & 1 & \text{} & \text{} & \text{} &
   \text{} & \text{} & \text{} & \text{} & \text{} &
   \text{} & \text{} & \text{} & \text{} &&&&&&&&&&& \\
 26 & 2 & 3 & 1 & 3 & \mathbf{{\underline{3}}}& \text{} & \text{} & \text{} &
   \text{} & \text{} & \text{} & \text{} & \text{} &
   \text{}  &&&&&&&&&&&\\
 27 & \text{} & 1 & \text{} & \text{} & \text{} &
   \text{} & \text{} & \text{} & \text{} & \text{} &
   \text{} & \text{} & \text{} & \text{} &&&&&&&&&&& \\
 28 & \text{} & 2 & 2 & \text{} & \text{} & \text{} &
   \text{} & \text{} & \text{} & \text{} & \text{} &
   \text{} & \text{} & \text{} &&&&&&&&&&& \\
 29 & \text{} & \text{} & \text{} & \text{} & \text{} &
   \text{} & \text{} & \text{} & \text{} & \text{} &
   \text{} & \text{} & \text{} & \text{} &&&&&&&&&&& \\
 30 & 1 & 2 & 2 & 6 & 6 & 6 & 6 & 6 & 6 & 2 & 2 & 2 &
   \text{} & \text{} &&&&&&&&&&& \\
 \hline
 \mbox{$\#'_k$} &  & & 312 & 142 & 67 & 44 & 29 & 28 & 16 & 15 & 13 & 7 & 6 & 5 & 5 & 4 & 4 & 3 & 3 & 2 & 2 & 1 & 1 & 1 & 1 \\
    N'_k &  & & 980 & 690 & 330 & 114 & 60 & 150 & 30 & 30 & 30 & 30 & 10 & 6 & 6 & 6 & 6 & 6 & 6 & 6 & 6 & 6 & 6 & 2 & 6 \\
 C_k &  & &\!\! \!\! 1000 & 700 & 450 & 300 & 200 & 200 & 150 & 150 & 150 & 150 & 100 & 100 & 100& 50 & 50 & 50 & 50 & 30 & 30 & 30 & 30 & 30 &30 \\
 \hline
\end{array}
}
\]
}
\caption{\label{maintable} The number $|Q_k(N)|$ of twist-classes of non-CM modular forms with 
rational coefficients, weight $k$, and minimal level $N$.  The number 
$|Q^u_k(N)|$ of such classes where rationality is unforced 
is 0 (regular type), $1$ (boldface), or $2$ (underlined boldface).}
\end{table}

   The computer algebra system {\em Magma} \cite{magma} lets one easily identify $Q_k(N)$ for $kN$ sufficiently small.  
 Table~\ref{maintable} presents the sizes $|Q_k(N)|$ for $k \leq 50$ and $N \leq 30$, as well as 
 related information.  We expect in particular that for $k \geq 18$, the
 table accounts for everything: 
 \begin{Conjecture} \label{mc}  For $6 \leq k \leq 50$ there is a largest $N_k$ with $|Q_k(N_k)| \geq 1$.  
 For $18 \leq k \leq 50$, this $N_k$ is either 30, 10, 6, or 2, as reported on Table~\ref{maintable}.  
  For $k \geq 52$, there are no non-CM newforms
 with rational coefficients at all.  
 \end{Conjecture} 
 \noindent Restricted to $N=1$, our conjecture is a weak version of the well-known Maeda conjecture 
 \cite{HM97}.    For $N>1$, our conjecture is similarly related to a natural generalization of the Maeda conjecture \cite{Tsa14,DT}.   The boldface entries on Table~\ref{maintable} reflect
  factorizations of Hecke polynomials that this generalization says never happens for $k$ sufficiently
  large.  Thus a novelty of  Conjecture~\ref{mc} is
 its effectivity: it says that for $k \geq 18$, there are only the indicated 
 exceptions for $N=3$, $15$, and $22$.  
 
 \subsection{Further discussion of the conjecture} Let $\#_k = \sum_{N=1}^\infty |Q_k(N)|$, which {\em a priori} 
 may be infinite or finite.   The situation changes as $k$ increases as follows.  
 By the modularity result of Wiles et al.\ \cite{Wil95,TW95,BCDT01}, the set
 $Q_2(N)$ is naturally identified with the set of twist-and-isogeny classes of 
 non-CM elliptic curves of minimal conductor $N$.  Via this connection, Cremona has identified the sets $Q_2(N)$
  for $N \leq 400,000$ \cite{LMFDB} and it easy to see that $\#_2$ is infinite.  The 
 case $k=4$ is similarly related to rigid
 Calabi-Yau three-folds \cite{GY11} and examples have been systematically pursued \cite{Mey05}; it seems to us 
 premature to speculate whether 
 $\#_4$ is infinite or finite.  
 The cases $k \geq 6$ have been studied \cite{Yui13,PR15}, but there do not seem to
 be any systematic non-modular sources: this lack of sources contributes to 
 our expectation of finiteness for these $\#_k$.
 
     The last block of Table~\ref{maintable} includes quantities 
  $N_k'$, and $C_k$.  Direct calculation is supportive of Conjecture~\ref{mc} as follows.   
 For $6 \leq k \leq 50$, we have computed $Q_k(N)$ for all $N \leq C_k$,
 finding it last non-empty at $N=N_k'$.   
 In weights $k=6$ and $8$,   
 non-empty $Q_k(N)$ become increasingly rare as $N$ approaches $C_k$, as illustrated
 by Figure~\ref{summatory}.    The thinning is rapid enough that
 we expect overall finiteness, although we also expect 
 that the observed maximum $N_k'$ may be considerably less than the conjectured actual maximum $N_k$.   
 In each of the weights $k=10$, $12$, $14$, and $16$, we think it more likely than
 not that $N_k'=N_k$.  In weights $18 \leq k \leq 50$, the ratio $C_k/N_k'$ is
 always at least five, giving us considerable confidence in $N_k'=N_k$, as asserted
 by the conjecture.  Similarly, we have carried our computations with cutoffs $C_k \geq 6$
 for $52 \leq k \leq 100$,  always finding $Q_k(N)$ to be empty.   
 
      The last block of the table also contains the lower bound 
 $\#_k' = \sum_{N=1}^{C_k} Q_k(N)$ to $\#_k$.  
Reformulating some of the previous discussion,
 our expectation is that $(\#_2,\#_4,\#_6,\dots)$ takes the form
 $(\infty, \#_4, 312'',  142'', 67', 44', 29', 28', 16, 15, \dots)$.   Here $\#_4$, whether
 it is infinite or finite, should be substantially large than $\#_6$.  Also $A''$ indicates
a number slightly larger than $A$, and $B'$ indicates a number either equal to 
or very slightly larger than $B$.    It seems likely that the sequence $\#_k$ is 
weakly decreasing, which would conform to general expectations in
an unexpectedly sharp way.

 \subsection{Content of the sections}     Section~\ref{review} gives the promised review of modular forms, using the case of $N=1$ and
  the familiar ring $M(1) = \C[E_4,E_6]$ as an example.  Included in 
  this section is a brief summary of the classification of CM newforms,
  which compares interestingly with our Conjecture~\ref{mc}.    Section~\ref{versus} discusses
 a decomposition 
 of $Q_k(N)$ into its ``forced'' and ``unforced" parts:
 \begin{equation}
 Q_k(N) = Q^f_k(N) \coprod Q^u_k(N).
 \end{equation}  
 It explains how Conjecture~\ref{mc} with $Q_k(N)$ replaced by $Q^f_k(N)$
 would become provable, with finiteness also for  $k=2$ and 
 $k=4$, and all  $N_k$ identifiable.  As to $Q_k^u(N)$, it
 gives a quantitative model,  supporting our expectation
 of emptiness for large $k$.  
  Section~\ref{assembling} explains the 
 simple calculations underlying Table~\ref{maintable}.   
  
 There is a large literature devoted to the explicit study of rational newforms 
 in weight two, mostly through their connection with elliptic curves.  
 Our viewpoint is that rational newforms in weight greater than two are worthy
 of at least a modest fraction of this detailed attention.  Sections~\ref{rings}
 and \ref{reps} are in this spirit, and study the cases
  $N=2$, $3$, $4$, $6$, and $8$.   Section~\ref{rings} focuses on
   rings $M(N)$ of modular forms, describing
  them via interesting overrings $\widehat{M}(N)$, 
  all of which are free on two generators like 
  $\widehat{M}(1) = M(1)$. Section~\ref{reps}
  presents many congruences between 
  newforms for a given $N$ and interprets these congruences in terms of explicit number fields.  
 Finally, Section~\ref{discussion} returns the focus to Conjecture~\ref{mc}, and 
 briefly discusses possible future
 directions.

\subsection{Acknowledgements.} The author thanks the conference organizers
for the opportunity to speak at {\em Automorphic forms: theory and computation} at King's College
  London, in September 2016.  This paper grew out of the first half of the author's 
  talk.  The list of newforms drawn up here is applied in \cite{RobPGL2}, which is 
  an expanded version of the second half.  The author's research was supported by grant \#209472 from
  the Simons Foundation and grant DMS-1601350 from the National Science Foundation.

\section{Review of modular forms}
\label{review}  This section presents a brief synopsis of the theory of modular forms,
as presented in more detail in e.g.\ \cite{Kob93} and \cite{St07}.     
Our purpose is to make this paper immediately accessible to a broad range of readers.     
We restrict to the case of trivial character until the very last subsection.  Throughout,
the classical case of level $N=1$ is used as an example.  

\subsection{Rings of modular forms}
\label{rings1}        For any pair $(k,N)$ consisting of a non-negative even integer weight $k$ and a positive integer level $N$, 
one has a corresponding finite-dimensional
complex vector space $M_{k}(N)$ of modular forms. 
These modular forms are functions on the upper half plane $\{z \in \C : \mbox{Im}(z) > 0\}$ 
satisfying certain transformation laws which become less demanding as $N$ becomes
more divisible.  In particular, the functions can be expressed as power series in $q = e^{2 \pi i z}$  and 
for us it suffices to simply regard all the spaces $M_{k}(N)$ as subspaces of the ring $\C[[q]]$
of formal power series in $q$. 

The sum of all these $M_k(N)$ together forms a graded ring $M(N)$.   Each
space $M_{k}(N)$ contains a subspace $S_{k}(N)$ of cusp forms, and
these cusp forms together form an ideal $S(N)$ in the ring $M(N)$. 
The ring $M(1)$ was studied in 1916 by Ramanujan \cite{Ram16}.  
To make room for later subscripting conventions, we use Ramanujan's
notations for certain Eisenstein series, writing $Q=E_4$, $R=E_6$.  
The ring then takes the form  
 $M(1) = \C[Q,R]$ with 
\begin{align*}
Q & = 1 + 240 \sum_{n=1}^\infty \sigma_3(n) q^n \!\!\!\!\! && =  1 + 240 q + 2160 q^2 + 6720 q^3 + \cdots \!\!\! & \in M_4(1), \\
R & = 1 - 504 \sum_{n=1}^\infty \sigma_5(n) q^n \!\!\!\!\!&& = 1 - 504q - 16632q^2 - 122976q^3 - \cdots \!\!\!& \in M_6(1).
\end{align*}
Here the formulas refer to the usual sum of positive divisors, $\sigma_j(n) = \sum_{d|n} d^j$.  
The ideal of cusp forms has generator
\begin{align*}
\Delta & = \frac{Q^3-R^2}{1728} = q - 24 q^2 + 252 q^3 - 1472 q^4 + 4830 q^5 - 6048 q^6 - \cdots.
\end{align*}
Let $\eta = q^{1/24} \prod_{n=1}^\infty (1-q^n)$.  Then one has the remarkable alternative 
expression $\Delta = \eta^{24}$.

 \subsection{Operators and newforms}  
 \label{operators}
  There are many important operators on the spaces $M_{k}(N)$.  Among these 
  are the push-up operators, corresponding to positive integers
 $t$.  The operator for $t$ takes the form $g = \sum a_n q^n$ 
 in $M_{k}(N)$ into the form $g_t = \sum a_n q^{tn} \in M_k(tN)$.   Also 
 playing an explicit role for us is the commuting family of Atkin-Lehner 
 involutions $w_{p^e}$ of the graded ring $M(N)$, one for each prime power $p^e$
 exactly dividing $N$.  Finally, one has a commuting family of Hecke operators $T_p$ on $S_k(N)$, indexed by
primes $p$ not dividing $N$.    These Hecke operators commute with the Atkin-Lehner operators and 
are given by the explicit formula  
\[
T_p(\sum a_n q^n) = \sum (a_{pn} + p^{k-1} a_{n/p}) q^n.
\]
Here $a_x$ is understood to be $0$ if $x$ is non-integral.  

A form $q+\cdots \in S_k(N)$ which is a basis for a one-dimensional
eigenspace of the Hecke operators is called a newform.  We let
$P_k(N) \subset S_k(N)$ be the set of newforms.  These newforms span the
new subspace $S^{\rm new}_k(N)$.
As $(M,g,t)$  runs over triples, with $M$ a divisor of $N$, $g$ a newform in $P_k(M)$, 
and $t$ a divisor of $N/M$, 
the push-ups $g_t$ form a basis of $S_k(N)$.  

Suppose $q=p^e$ exactly divides $N$ and $g = \sum a_n q^n$ is a newform in $P_k(N)$ on which $w_{q}$ acts
with the eigenvalue $\epsilon_q$.    If $e=1$, then one has the formula 
\begin{equation}
\label{alsign}
a_p = - \epsilon_p p^{k/2-1}.
\end{equation}
This formula is very useful because it identifies $\epsilon_p$.  If $e>1$, then one has the simpler but not so useful 
formula $a_p=0$.  

The ring $M_k(1)$ is different from all the other $M_k(N)$ in that no push-up or Atkin-Lehner operators
are involved in its description.   However the $T_p$ behave completely typically.  For all $N$ and 
$g = \sum_n a_n q^n \in P_k(N)$, one has the simple formula $T_pg = a_pg$.  Also if $m$ and $n$ 
are relatively prime then $a_{mn} = a_m a_n$.  For example, the $a_6 = -6048$ for $\Delta \in P_{12}(1)$
is indeed $a_2 a_3 = -24 \cdot 252$.
 
  \subsection{Dimension formulas}  We are most directly interested in the spaces $S^{\rm new}_k(N)$. 
There is an exact formula \cite[Prop~6.1]{St07} for the dimension of the larger space $S_k(N)$.
Taking the  largest term as an approximation gives
 \begin{equation}
 \label{fulldimension}
 \dim S_k(N) \approx \frac{k-1}{12} \prod_{p^e||N} p^{e-1} (p+1).
 \end{equation}
%
The exact general
 formula for $S_k(N)$ passes to one for $S_k^{\rm new}(N)$.  The approximate 
 formula becomes 
 \begin{equation}
 \label{newdim}
  \dim S^{\rm new}_k(N) \approx \frac{k-1}{12} \prod_{p^e || N} m(p,e),
 \end{equation}
  with
  \[
  m(p,e) = \left\{ \begin{array}{ll} 
   p-1 & \mbox{if $e=1$,} \\ p^2-p-1 & \mbox{if $e=2$,} \\
  (p-1)^2 (p+1) p^{e-3} & \mbox{if $e \geq 3$.} \end{array} \right.
  \]
  Because the spaces $S_k(d)$ are involved for all $d|N$, 
  secondary terms are more complicated in the exact formula for  $\dim S^{\rm new}_k(N)$.
  
    Fix $N = \prod_{i=1}^m q_i$ with $q_i = p_i^{e_i}$ and 
     $p_1 < \cdots < p_m$.  Then the
 Atkin-Lehner operators give decompositions
 \begin{align*}
 M_k(N) & = \sum_\epsilon M_k(N)^\epsilon, &
S_k(N) & = \sum_\epsilon S_k(N)^\epsilon, &
S^{\rm new}_k(N) & = \sum_\epsilon S^{\rm new}_k(N)^\epsilon. 
 \end{align*}
 Here $\epsilon$ runs over the $2^m$ sign strings $(\epsilon_{q_1},\dots,\epsilon_{q_m})$,
 with $w_{q_i}$ acting by $\epsilon_{q_i}$.   As one might expect, an approximate
 formula for $S_k(N)^\epsilon$ is $1/2^m$ times \eqref{fulldimension}.  
 However for $S^{\rm new}_k(N)^\epsilon$, one has to replace 
 each $m(p_i,e_i)$ by the appropriate $m^{\epsilon_{q_i}}(p_i,e_i)$.  
 Here
 \begin{equation}
\label{localmasses}  m^\pm(p,e) = \left\{ \begin{array}{ll} 
   (p-1)/2, & \mbox{if $e=1$,} \\ 
   (p^2-p -1 \mp 1)/2, & \mbox{if $e=2$,} \\
  (p-1)^2 (p+1) p^{e-3}/2,& \mbox{if $e \geq 3$.} \end{array} \right.
 \end{equation}
For fixed $N$ and increasing $k>2$, all approximations discussed in this section
are off by a function which is periodic in $k$.    An interesting feature of 
\eqref{localmasses} is $m^{+}(2,2)=0$, and indeed $P_k(N)^\epsilon$ is empty
whenever $\ord_2(N) = 2$ and $\epsilon_4=+$.  

As an example, from the description of $M(1)$ and $S(1)$ above, one has 
the exact formula
\[
\dim(S_k(1)) = \frac{k-1-\delta_k}{12},
\]
with $\delta_k = 13$, $3$, $5$, $7$, $9$, $-1$ for $k=2$, $4$, $6$, $8$, $10$, $12$
and satisfying $\delta_{k+12} = \delta_{k}$.  We are interested particularly
in one-dimensional spaces.  These occur for $k=12$, $16$, $18$, $20$, $22$, and $26$.
The corresponding unique newforms are
 \begin{align}
 \label{raman}
 \Delta, &&
 Q \Delta, &&
 R \Delta &&
 Q^2 \Delta && 
 QR \Delta, &&
 Q^2 R \Delta.
  \end{align}  
 These newforms were all studied by Ramanujan \cite{Ram16}; they provide six
 explicit illustrations of the objects of our title, {\em Newforms with rational 
 coefficients.}

\subsection{Quadratic twists}    
\label{quadtwists}  Quadratic number fields are classified by their discriminants 
$D$.  These integers, and also the discriminant $1$ of the quadratic algebra $\Q \times \Q$,
are called fundamental discriminants.  Explicitly, a fundamental discriminant is 
an integer of the form $td$, where $d$ is a square-free integer congruent to $1$ modulo $4$, 
and $t \in \{1,-4,8,-8\}$.  Fundamental discriminants form a set of representatives in 
$\Q^\times$ of the group $\Q^\times/\Q^{\times 2}$.    Each fundamental discrimininant
$D$ gives a character $\chi_D : (\Z/D)^\times \rightarrow \{-1,1\}$ given
by the quadratic residue symbol, $\chi_D(n) = (D/n)$.  

The infinite group of these quadratic Dirichlet characters 
acts on the set of newforms of a given 
weight $k$ by twisting.    Suppose $g = \sum_{n=1}^\infty a_n q^n \in P_k(N)$ 
and $\chi = \chi_D$.  Then $g^{\chi}$ is a newform with 
level $N_{g^\chi}$ 
dividing
$\LCM(N,D^2)$.  It is characterized by $a_n^\chi = \chi(n) a_n$ for $n$
not dividing $\LCM(N,D^2)$.   Also, equality holds in
$\ord_p(N_{g^\chi}) \leq \max(\ord_p(D^2),\ord_p(N))$ if
$\ord_p(D^2) \neq \ord_p(N)$.  

For a form $g$, one says its {\em minimal level} is the minimum of the
levels of all its twists.   Say that a form is {\em minimal} if
its level is equal to its minimal level.  If $g$ is minimal
of level $N$ then $g$ has $t(N)$ minimal twists, where
$t$ is the multiplicative function satisfying 
\begin{equation*}
 t(2^e)  = 
 \left\{
 \begin{array}{ll} 
 1 & \mbox{if $e \in \{0,1,2,3\},$} \!\!\!\!\! \\
 2 & \mbox{if $e \in \{4,5\}$,} \\
 4 & \mbox{if $e \geq 6$.} \\
 \end{array}
 \right.
 \mbox{ and, for $p$ odd, }
  t(p^e)  = 
 \left\{
 \begin{array}{ll} 
 1 & \mbox{if $e \in \{0,1\}$,} \\
 2 & \mbox{if $e \geq 2$.} 
 \end{array}
 \right.
\end{equation*}
 The naturality of this definition is apparent from its alternative description:
 $t(N)$ is the number of fundamental discriminants $D$
 such that $D^2|N$. As a variant of the standard notion of squarefree, say 
that an integer $N$ is {\em quadfree} if $t(N) = 1$.    So 
$N$ is quadfree if $e := \ord_2(N) \leq 3$ and 
its odd part $N/2^e$ is squarefree.   

\subsection{CM newforms and the set $Q_k(N)$ of interest}  
\label{CMnewforms} 
We can now define
some of the terms used in the introduction.    Most newforms $g$ satisfy
$g = g^\chi$ only for the trivial character $\chi$.   The remaining newforms
satisfy $g = g^\chi$ for exactly one non-trivial character $\chi = \chi_D$
and moreover $D$ is negative.  Such a newform is said to have CM by $D$.  
Necessarily, its level is divisible by $D^2$.     
So for most $N$, we do not encounter CM newforms; for a very 
few $N$ we do, and then we just discard them.   

In practice, the CM newforms to be discarded are immediately recognizable by their Fourier expansions:
while coefficients $a_p$ are rarely or perhaps even never zero for a non-CM newform, they
are always zero whenever $\chi_D(p)=-1$ for newforms with CM by $D$.  On 
a rigorous level, we can confirm that a newform with apparent CM by $D$ really
does have CM by $D$ by a general structure theorem \cite[Theorem~2.4]{Sch09}.  
Namely quadratic twist-classes of CM newforms 
in weight $k = 1+e$ are in bijection with
imaginary quadratic fields $\Q(\sqrt{D})$ with class group
of exponent dividing $e$.  Also minimal levels are known \cite[Table~1]{Sch09}.  

The classification just mentioned lets one see that the finiteness 
assertions of Conjecture~\ref{mc} are true to some 
extent in the parallel situation of CM newforms, but false 
in their full statement, as follows.  
For fixed $k$, the complete list of minimal levels $N$ of CM newforms is sometimes known.  For
example for $k=2$, the discriminants $D$ include $-3$, $-4$,  and
$-8$, with associated minimal levels $N = 27$, $32$ and $256$.  
The remaining discriminants are $D=-p$ with $p \in \{7,11,19,43,67,163\}$,
always with minimal level $N=p^2$.   The twenty-six
possible $D$ for $k=4$ and their levels $N$ are likewise listed in \cite[Table~3]{Sch09}.  
For general fixed $k$, finiteness of the list of $N$ is expected,
as it is implied by the Riemann hypothesis for $L$-functions
of odd real Dirichlet characters \cite[Theorem~2.1]{Sch09}.
For fixed $N$ on the other hand, the set of weights $k$ for which there
is at least one rational CM newform with minimal level $N$ can easily be infinite.  
For example, this set of
weights is all positive even integers for the six $N=p^2$ above.

Returning to our main focus, all terms in the definition of $Q_k(N)$ from the introduction 
have now been defined: $Q_k(N)$ is the set of all twist-classes of
non-CM newforms with minimal level $N$ and rational coefficients. 
Such a class $[g]$ has exactly $t(N)$ representing newforms of level $N$.
In the common case that $N$ is quadfree, it has just one
such representative.  

\subsection{General character}  
\label{genchar}
It is standard to work in greater
generality, defining spaces $M_k(N,\chi)$ for general Dirichlet 
characters $\chi$ and general weights $k \in \Z_{\geq 0}$.  
These spaces can be nonzero only if 
the conductor of $\chi$ divides $N$ and $\chi(-1) = (-1)^k$.
We have summarized standard material in the 
setting of trivial character only, meaning our
$M_k(N)$ agrees with $M_k(N,\chi_1)$.   

We are focusing on the case of trivial character because
a non-CM newform with rational coefficients necessarily has
trivial character.    In Section~\ref{rings} below, we will encounter cases
of non-trivial character.  We will even allow weights
to become half-integral.   However our
summary here is adequate for following the discussion
there.

\section{Forced vs.\ unforced rationality}
\label{versus}
    This section discusses forced versus unforced rationality. 
We include a small {\em Magma} program in this section and two
more in the next.    These programs allow readers not familiar with
{\em Magma} to directly see some of the calculations
 underlying this paper.   Run on small enough arguments,
the programs finish in less than the two minutes allowed
by the free online {\em Magma} calculator.  

\subsection{The Galois action}  
\label{Galoisaction}
Fix a weight $k \in 2 \Z_{\geq 1}$ and a level $N \in \Z_{\geq 1}$.  
The group $\Gal(\overline{\Q}/\Q)$ acts on the corresponding set $P_k(N)$ of newforms
by conjugating coefficients.   Let $p$ be a prime number not dividing
 $N$.  Let $f_{k,N,p}(x)$ be the characteristic 
polynomial of $T_p$ acting on the space $S_k(N)^{\rm new}$.  
Then, assuming $f_{k,N,p}(x)$ is separable, the action
of $\Gal(\overline{\Q}/\Q)$ on the set $P_k(N)$ agrees
with its action on the coefficients $a_n$, these being the roots
of $f_{k,N,p}(x)$.

A general program computing the $f_{k,N,p}(x)$ is obtained by
concatenating built-in {\em Magma} commands.  First, to obtain
output in a standard form, one can introduce the variable
$x$ by \verb@_<x>:=PolynomialRing(Integers());@.  The general program is then
\smallskip

\verb@charpol := func<k,N,p|@

\verb@Factorization(CharacteristicPolynomial(@

\verb@HeckeOperator(NewSubspace(CuspForms(N,k)),p)@

\verb@))>;@
\smallskip

\noindent To compute say $f_{50,3,2}(x)$, one inputs all of the above and then
\verb@charpol(50,3,2);@.  In about a second, one finds that $f_{50,3,2}(x)$ 
factors as an irreducible quartic times an irreducible quintic.  

    The action of $\Gal(\overline{\Q}/\Q)$ on $P_k(N)$ passes to an action on $Q_k(N)$.
In this paper we are interested in the fixed points of $\Gal(\overline{\Q}/\Q)$ on $Q_k(N)$.  
When $N$ is quadfree, as defined at the end of \S\ref{quadtwists}, it is simplest to think of $Q_k(N)$ 
as being simply its single representing form in $P_k(N)$.  
Similarly for general $N$, a fixed point on $Q_k(N)$ can be
thought of using \S\ref{quadtwists} as $t(N)$ fixed points, 
all twists of one another, on $P_k(N)$.   

    There are different ways of expressing the Galois action.  One convenient
way is to let $E_{k,N} \subseteq \mbox{End}(S^{\rm new}_k(N))$ be the 
$\Q$-algebra generated by all the Hecke operators $T_p$.  One has $E_{k,N} = \Q[x]/f_{k,N,p}(x)$
whenever $f_{k,N,p}(x)$ is separable.  In this case the factorization
of $E_{k,N}$ into fields, which is the issue under study, is exactly reflected
in the factorization of $f_{k,N,p}(x)$ into irreducible polynomials.  

\subsection{The case $N=1$}   The Maeda conjecture \cite[Conj~1.2]{HM97} says that 
the image of $\Gal(\overline{\Q}/\Q)$ in its action on $P_k(1)$ is always
the full symmetric group on the degree $d_k$.  A slightly strengthened version 
\cite[Conj~1.1]{GM12} includes also the separability of all $f_{k,1,p}(x)$.  In other
words, it says that the Galois group of $f_{k,1,p}(x)$ is always 
$S_{d_k}$.  The reference \cite{GM12} also proves the strengthened conjecture for $p=2$ and 
$k \leq 14000$, and surveys other results related to the conjecture.  

In the six cases when $d_k=1$, the set $P_k(1)$ coincides with its
subset of rational forms, these forms having been listed in \eqref{raman}.  For the 
$k$ with $d_k \geq 3$, the Maeda conjecture is a stronger statement
than Conjecture~\ref{mc}'s assertion that $Q_k(1)$ is empty.  The computed
 Hecke fields $E_{k,1} = \Q[x]/f_{k,1,2}(x)$ have 
very large discriminants, rapidly increasing with $k$.  For example, the first case beyond
rationality is the quadratic field $E_{24,1}$, and its discriminant is already the fairly large
 prime number $144169$.   The large discriminants constitute further heuristic evidence 
 that $f_{k,1,2}(x)$ is always irreducible.

\subsection{Quadfree $N$}    For general $N$, the 
action of $\Gal(\overline{\Q}/\Q)$ on $P_k(N)$ stabilizes
the Atkin-Lehner subsets $P_k(N)^\epsilon$.   
A very plausible analog of the Maeda conjecture 
for general $N$ was formulated in \cite{Tsa14} and strengthened 
to be numerically more precise and include Galois groups 
in \cite{DT}.   For the case of quadfree $N$,
 it essentially says that, a finite number of 
exceptional spaces aside, each $P_k(N)^\epsilon$
behaves qualitatively like $P_k(1)$.   More precisely,
let $p$ be the smallest prime not dividing $N$.  Then, outside of
finitely many $k$,  the characteristic polynomial $f_{k,N,p}^\epsilon(x)$ should
have Galois group the full symmetric group on its degree.   

As an example of generic behavior, consider again the polynomial
$f_{50,3,2}(x)$ computed in \S\ref{Galoisaction}.   The Galois groups
of its irreducible factors 
are $S_4$ and $S_5$, and the field discriminants have 51 and 
79 digits.  The factorization is completely expected as 
the summands in the decomposition $P_{50}(3) = P_{50}(3)^+ \coprod P_{50}(3)^-$ 
have size $4$ and $5$ respectively.  

\subsection{General $N$} \label{generalN}  For $N$ which are not quadfree,
 there are structures on the set $P_k(N)$ which go beyond the
decomposition induced by the Atkin-Lehner operators.
To give an indication of this phenomenon, 
consider the case $N=25$.     Then $P_k(25)$ breaks
into six parts, which are conveniently described 
via twisting by $\chi = \chi_5$.    Write 
$g \in P_{k}(25)^{\renewcommand{\arraystretch}{.08} \begin{array}{c} \!\!\! {}_\epsilon \\ \; \\ \! \!\! {}_\omega \end{array}} \!\!\!$ if $g \in P_k(25)^\epsilon$ and 
$g^\chi \in P_k(25)^\omega$.  Then
\[
P_k(25) = P_{k}(25)^{\pm} \coprod  P_k(25)^{\mp} \coprod P_k(25)^{=} \coprod \mbox{Rest}.
\]
Here we do not particularly care about $\mbox{Rest}$, as it consists of the
forms $g \in P_k(25)^+$ with $f^\chi \in P_k(1)$, $P_k(5)^+$, or 
$P_k(5)^-$.    The number of possibilities for $N=25$ replaced by 
a general $p^e$ is given in \cite[Prop 3]{DT}.

In general, one has a decomposition of $P_k(N)$ into $\Gal(\overline{\Q}/\Q)$ stable summands $P_k(N)^\delta$ refining
the Atkin-Lehner involution.   Just as $\epsilon$ is a string of signs indexed by the prime powers exactly dividing $N$,
so too is $\delta$ a string of symbols indexed by these prime powers.  For example, take $N=150 = 2 \cdot 3 \cdot 5^2$.
The set $P_k(150)$  breaks into $2 \cdot 2 \cdot 6 = 24$ parts, with $2 \cdot 2 \cdot 3=12$ corresponding
to minimal level $150$.    For large $k$,  one can expect $f_{k,150,7}(x)$ to factor
into twenty-four irreducible polynomials, one for each type $\delta$.  However $f_{16,150,7}(x)$ factors
into twenty-five irreducible polynomials.  Here the number of irreducible polynomials arises as
$25=24-1+2$, as follows.  Not so interestingly,  $P_{16}(6)^{+-}$ being empty implies that one of the $\delta$ is 
not represented.  However there is a rational newform with type $\delta' = (+,+,\mp)$, namely
\begin{equation*}
g =  q - 128 q^2 - 2187 q^3 + 16384 q^4 + 279936 q^6 - 511994 q^7 + \cdots.
\end{equation*}
Its twist $g^{\chi_5}$ has type $\delta''=(-,-,\pm)$.   These two exceptional $\delta$ each contribute 
an extra irreducible factor, as   
 \begin{eqnarray}
\nonumber \lefteqn{ f_{16,150,7}^{+,+,\mp}(x) = (x + 511994) \cdot } \\
\label{threeplusone} &&(x^3 + 701247x^2 - 5978366987397x + 3322646963771081149)
 \end{eqnarray}
 and $f_{16,150,7}^{-,-,\pm}(x) =  f_{16,150,7}^{+,+,\mp}(-x)$.  
 
 Many $P_k(N)^\delta$ cannot possibly have $\Gal(\overline{\Q}/\Q)$ fixed points, because the structures on $P_k(N)^\delta$
 imply that all factors of the Hecke algebra $E_{k,N}^\delta$ contain a specified cyclotomic field 
 larger than $\Q$.  This phenomenon is familiar from the case of weight $2$ forms and elliptic curves.  It holds
 without change for $k \geq 4$.  In particular, rationality of a form in $P_k(N)^\delta$ implies that $\ord_2(N) \leq 8$,
 $\ord_3(N) \leq 5$, and $\ord_p(N) \leq 2$ for larger primes $p$.  
 
 \subsection{The dichotomy}
 \label{dichotomy}
 We can now explain the decomposition of $Q_k(N)$ into its two parts $Q^f_k(N)$ and $Q^u_k(N)$. 
 Let $[g] \in Q_k(N)$ with representing newform $g \in P_k(N)$.  We say that the rationality
 of $g$ is {\em locally forced}, or simply {\em forced},
  if $g$ is the only non-CM newform in its refined part $P_k(N)^\delta$.   
 We write $[g] \in Q^f_k(N)$ in this case and $[g] \in Q^u_k(N)$ otherwise.   
 Exact formulas for $|P_k(N)^\delta|$ are not available at the moment, but they are surely within
 reach.  Using these formulas and also formulas or upper bounds 
 for the number of CM newforms, one could
 compute a cutoff $c_k$ such that $Q^f_k(c_k)$ is nonempty but $Q^f_k(N)$ is empty for $N > c_k$.    
 In fact, we are asserting in Conjecture~\ref{mc} that for $k \in [18,50]$, the 
 cutoff $c_k$ is the number listed as $N_k'$ on Table~\ref{maintable}.   
 Similarly, the next paragraph gives strong evidence that $c_{16}=42$, and 
 \eqref{sequenceends} below suggests further than $(c_{14},c_{12},c_{10}) = (42, 90, 210)$.

 To see 5=2+2+1 unforced instances of rationality, consider $k=16$.  
 As reported on Table~\ref{maintable}, rationality is forced for all the newforms with $N \leq 30$,
 except for the two with $N=8$.  For just one $N \in [31,149]$ is the set  $Q_{16}(N)$ non-empty,
 namely $N=42$, where it has size four.   For this exceptional level,   
 $f_{16,42,5}(x)$ factors into four linear and five quadratic irreducible
 factors.   Closer inspection shows that $|P_{16}(42)^\epsilon|=1$ for $\epsilon = (-,-,-)$ 
 and $(-,+,+)$, while otherwise $|P_{16}(42)^\epsilon|=2$.  Thus $|Q^f_{16}(42)| = |Q^u_{16}(42)|=2$.
 In fact, the source of the unforced rationality is 
 \[
 f_{16,42,5}^{+,+,-}(x) = (x+58290)(x-296442).  
 \]
  For the example presented in the last subsection,
 the unexpected factorization \eqref{threeplusone} is saying that 
 $|Q^f_{16}(150)|=0$ and 
 $|Q^u_{16}(150)|=1$.

\subsection{A heuristic model}
This subsection describes a heuristic model 
for the factorization of the Hecke algebras $E_{k,N}^\delta$,
obtained by considering the factorization of defining polynomials 
$f^\delta_{k,N,p}(x)$.
While the model is very rough, we feel it complements 
our catalog of newforms by supporting Conjecture~\ref{mc} in
a different way.

Fix $(k,N,\delta)$ and a prime $p$ not dividing $N$.   
The quantity $w = p^{(k-1)/2}$ plays the role of a scaling factor.  
Let $d = |P_k(N)^\delta|$, so that the monic polynomial
$f^\delta_{k,N,p}(x) \in \Z[x]$ has degree $d$.   Its
$d$ roots are all real with absolute value at most
$2w$.  The approximate number of such polynomials
has the form 
\begin{equation}
\label{volume}
V_d(w) = \left[\frac{2^d}{d!} \prod_{j=1}^d \left(\frac{2j}{2j-1}\right)^{d+1-j} \right]  w^{\Delta(d)}.
\end{equation}
Here $\Delta(d) = d(d+1)/2$ is the $d^{\rm th}$ triangular number and 
the complicated coefficient is a volume computed in \cite[Prop~2.2.1]{DH98}.

Now let $r<s$ be positive integers summing to $d$.   The chance that a polynomial
in the ensemble of degree $d$ polynomials under consideration factors 
into a degree $r$ polynomial times a degree $s$ polynomial is approximately 
$\mbox{Prob}_{r,s}(w) = V_r V_s/V_d$.  If $r=s$, then this formula double counts, and
one needs to insert a $2$ in the denominator.  Applying \eqref{volume} three times
and simplifying, 
one gets 
\begin{equation}
\mbox{Prob}_{r,s}(w) = \left[ \frac{d!}{2^{\delta_{rs}} r! s!} \prod_{j=1}^d \left( \frac{2j-1}{2j} \right)^{j-1} \right] \frac{1}{w^{rs}}.
\end{equation}
For example, the first non-trivial case is $r=s=1$.  Here the chance that a
quadratic polynomial from the ensemble splits is approximately $3/(4w) = 3/(4p^{(k-1)/2})$.   

\begin{figure}[htb]
\begin{center}
\includegraphics[width=4.7in]{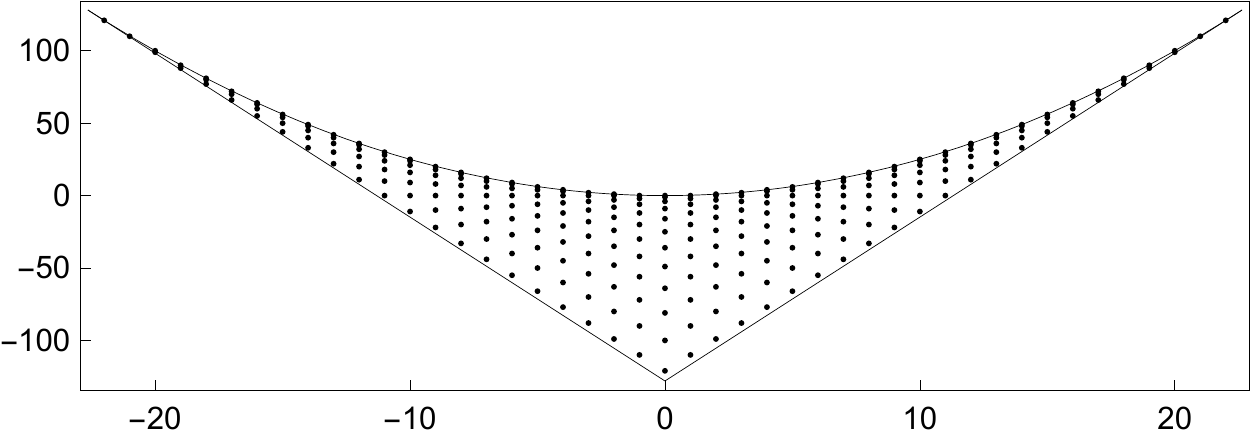}
\end{center}
\caption{\label{shield}  Illustration in the $b$-$c$ plane of the heuristic that a quadratic 
Hecke polynomial for the prime $p=2$ and weight $k=6$ has
approximately approximately a $1/8$ chance of factoring.  
}
\end{figure}

The first non-trivial case in more detail goes as follows.  
 The factorizing polynomials are simply $(x-a_1)(x-a_2)$ with $|a_i| \leq 2w$
and $a_1 \leq a_2$. All allowed polynomials are $x^2+bx+c$ with $|b| \leq 4 w$ and 
$w(|b|-w) \leq c \leq w^2/2$.   The regions of the $a_1$-$a_2$ and  $b$-$c$ planes 
given by the inequalities respectively have area $V_1(w)^2/2 = 8w^2$ and  $V_2(w) = 32 w^3/3$.  Their 
ratio is indeed $P_{1,1}(w) = 3/(4 w)$.   Figure~\ref{shield}
draws this case for $(p,k)=(2,6)$ so that $w = 2^{5/2}$.  Factoring polynomials
are represented as points.  There are
$276$ of them, which is indeed close to the approximation
 $V_1(w)^2/2 = 256$.  The total number of 
polynomials $x^2+bx+c$ is $1951$ while $V_2(w)  \approx 1930.9$.   

One way that our model is very rough is that it does not
account for the fact that roots of Hecke polynomials should be distributed
according an approximation of the Sato-Tate measure.  More
seriously, the number $\mbox{Prob}_{r,s}(p^{(k-1)/2})$ depends on 
$p$, whereas the factorization behavior of $E^\delta_{k,N}$
is independent of $p$.   To proceed further, we take as our
heuristic that the chance of $E^\delta_{k,N}$ factoring
into a degree $r$ algebra times a degree $s$ algebra is
$\mbox{Pr}_{r,s}(k) := \mbox{Prob}_{r,s}(2^{(k-1)/2})$, 
whether or not $2$ divides $N$. 

Some explicit numerics are as follows.   Corresponding to three of the $\underline{\bf 2}$'s on Table~\ref{maintable},
factorization patterns for $E^+_{6,17}$, 
$E^-_{16,8}$, and $E^+_{22,3}$ are all $1+1$.  Corresponding splitting probabilities are $\mbox{Pr}_{1,1}(k)$
for $k=6$, $16$, and $22$, namely $13.3\%$, $0.4\%$, and $0.05\%$.   The other $\underline{\bf 2}$ on 
Table~\ref{maintable} comes from $E^+_{6,19}$ and $E^{-}_{6,19}$ having splitting behavior $2+1$ 
and $4+1$ respectively.  The probability here is 
$\mbox{Pr}_{2,1}(k) \mbox{Pr}_{4,1}(k) \approx 5.8\% \cdot 0.24\% = 0.014\%$.    

Some qualitative phenomena which the heuristic seems to get right are as follows.  
First, the chance of splitting decreases rapidly as the weight $k$ increases.  Second, for fixed $k$, this chance 
also decreases rapidly with increasing degree $d$.   Third, factorizations of the form $d = 1 + (d-1)$ 
are much more common than factorizations of any other type.   In fact, 
we have observed no other factorizations for quadfree $N$ in weights 
$k \geq 6$  except for $|P_6(3 \cdot 23)^{+-}| = 2+4$, $|P_6(3 \cdot 5 \cdot 17)^{++-}| = 3+4$, $|P_6(455)^{+--}| = 2+10$,  and $|P_8(3 \cdot 17)^{--}|=2+4$. 
This third phenomenon is one of the reasons that this paper concentrates not on general
factorizations, but rather on just those which produce newforms with 
rational coefficients. 

On a more quantitative level, there are definitely more factorizations
than the heuristic predicts.  In weight two, the heuristic correctly
predicts that there are infinitely many elliptic curves, but 
considerably underestimates the number of elliptic curves
per level.   As examples in higher weight, 
one of the rational newforms discussed in  \S\ref{SampleLarger} 
has associated probability $\mbox{Pr}_{1,12}(10) \approx 2.2 \times 10^{-16}$
while one from \S\ref{extravan} corresponds to the even smaller number 
$\mbox{Pr}_{1,83}(4) \approx 3.4 \times 10^{-37}$.  
One is thus led to ask for a conceptual source of 
rational newforms such as these, something
we will briefly pursue in \S\ref{extravan}.

\section{Assembling the main table}  
\label{assembling} In this section we discuss how we drew up 
 Table~\ref{maintable}, as well as its unprinted
extension to a larger region in the $N$-$k$ plane.   
We present a number of examples similar to those 
of the last subsection, but now with
more reference to Table~\ref{maintable} itself.  

\subsection{Computing the cardinality $|Q_k(N)|$}  
\label{cardinality} Combining several built-in {\em Magma} functions,
we define a new one:
\smallskip

\verb@RatNewforms := func<k,n|@

\verb@[fs : fs in Newforms(CuspForms(Gamma0(n),k))|#fs eq 1]>;@
\smallskip

\noindent Following this definition by \verb@RatNewforms(22,3)@ then returns approximations to the two elements in $Q_{22}(3)$:
\begin{eqnarray*}
f_1 & = &  q + 1728 q^2 - 59049 q^3 + 888832 q^4 - 41512770 q^5 - \cdots, \\
f_2 & = &  q - 2844 q^2 - 59049 q^3 + 5991184 q^4 + 3109950 q^5 + \cdots.
\end{eqnarray*}
This particular computation takes less than a second and 
accounts for the $\underline{\bf{2}}$ in the row $3$ column $22$ of Table~\ref{maintable}.    On the other hand, for $N$ a highly factorizing
level close to our cutoff $C_k$, \verb@RatNewforms(@$k$\verb@,@$N$\verb@)@ takes around an hour.  

To convert the output of \verb@RatNewforms(@$k$\verb@,@$N$\verb@)@ into cardinalities $Q_k(N)$
one needs to take quadratic twisting into account.  If no CM newform is present, for example if
a prime exactly divides $N$, then this cardinality is the number of newforms
returned, divided by the number $t(N)$ of twists allowed from \S\ref{quadtwists}.  When CM newforms are present, 
they are in practice easily recognized by their $a_p$ being zero for many $p$.  We
used the theory of CM newforms as summarized in \S\ref{CMnewforms}
 to confirm that CM is really present, discarded the forms,
and then divided the number of remaining forms by $t(N)$.

\subsection{Identifying the decomposition $Q_k(N) = Q^f_k(N) \coprod Q^u_k(N)$.}  In the example of the previous subsection,
 both coefficients of $q^3$ are $-3^{10}$, so by \eqref{alsign} they both belong to $S^{\rm new}_{22}(N)^+$.  They thus
 contribute to the unforced part $Q^u_{22}(3)$ of $Q_{22}(3)$.   To systematically separate the forced part from the unforced 
part in the case of quadfree $N$, we used the following refinement of our previous \verb@charpol@: 
\smallskip

\verb@charpol2 := function(k,N,signs,p)@

\verb@fN := Factorization(N); fN2 := [g[1]^g[2]:g in fN];@

\verb@New := NewSubspace(CuspForms(N,k));@

\verb@return Factorization(CharacteristicPolynomial(@

\verb@HeckeOperator(New,p)*(&*[(AtkinLehnerOperator(New,fN2[i])^2+@

\verb@signs[i]*AtkinLehnerOperator(New,fN2[i]))/2: i in [1..#fN2]])));@

\verb@end function;@

\smallskip

\noindent Running  \verb@charpol2(22,3,[1],2)@ and then \verb@charpol2(22,3,[-1],2)@ lets one conclude
$\det(x-T_2|S_{22}(3)^+) = (x - 1728)(x+2844)$ and $\det(x-T_2|S_{22}(3)^-) = x^2 - 666 x - 2464992$.
This computation shows again that rationality of these two newforms is unforced.    
Some non-quadfree cases can also be done via this program, as illustrated by the 
example of $(k,N) = (16,150)$ in \S\ref{generalN}.   

\subsection{Sample calculations at small level}  
\label{SampleSmall}
Table~\ref{dimensions} gives the result of running  \verb@charpol2@ 
for three hundred different $(k,N,\epsilon)$.   All polynomials obtained  were
irreducible, except for the cases $(22,3,+)$ and $(16,8,-)$ where
  \begin{table}[htb]
{
\[
{\renewcommand{\arraycolsep}{1.3pt}
\begin{array}{|rcc|c|ccccccccccccccccccccccccc|}
\hline
N&\epsilon_q&\epsilon_3&m&  2 & 4 & 6 & 8 & 10 & 12 & 14 & 16 & 18 & 20 & 22 & 24 & 26 & 28 & 30 & 32 & 34 & 36 & 38 & 40 & 42 & 44 & 46 & 48 & 50 \\   
\hline
1 &&&1&   &  &  &  &  & \bf1 &  &\bf1 & \bf1 & \bf1 & \bf1 & \color{g}{2} & \bf1 & \color{g}{2} & \color{g}{2} & \color{g}{2} & \color{g}{2} & \color{g}{3} & \color{g}{2} & \color{g}{3} & \color{g}{3} & \color{g}{3} & \color{g}{3} & \color{g}{4} & \color{g}{3} \\
 \hline
2&+&&0.5 &  &  &  & \bf1 & &  & \bf1& \bf1&  & \bf1& \bf1& \bf1& \bf1& \bf1& \bf1& \color{g}{2} & \bf1& \bf1& \color{g}{2} & \color{g}{2} & \bf1& \color{g}{2} & \color{g}{2} & \color{g}{2} & \color{g}{2} \\
2&-&&0.5&   &  &  &  & \bf1&  & \bf1&  & \bf1& \bf1& \bf1&  & \color{g}{2} & \bf1& \bf1& \bf1& \color{g}{2} & \bf1& \color{g}{2} & \bf1& \color{g}{2} & \color{g}{2} & \color{g}{2} & \bf1& \color{g}{3} \\
 \hline
3&&+&1 &   &  &  & \bf1& \bf1& \bf1& \bf1& \bf1& \bf1& \color{g}{2} & \bf2 & \color{g}{2} & \color{g}{2} & \color{g}{2} & \color{g}{2} & \color{g}{3} & \color{g}{3} & \color{g}{3} & \color{g}{3} & \color{g}{3} & \color{g}{3} & \color{g}{4} & \color{g}{4} & \color{g}{4} & \color{g}{4} \\
3&&-& 1&   &  & \bf1&  & \bf1&  & \color{g}{2} & \bf1& \color{g}{2} & \bf1& \color{g}{2} & \bf1& \color{g}{3} & \color{g}{2} & \color{g}{3} & \color{g}{2} & \color{g}{3} & \color{g}{2} & \color{g}{4} & \color{g}{3} & \color{g}{4} & \color{g}{3} & \color{g}{4} & \color{g}{3} & \color{g}{5} \\
 \hline
4&-&&1&   &  & \bf1&  & \bf1& \bf1& \bf1& \bf1& \color{g}{2} & \bf1& \color{g}{2} & \color{g}{2} & \color{g}{2} & \color{g}{2} & \color{g}{3} & \color{g}{2} & \color{g}{3} & \color{g}{3} & \color{g}{3} & \color{g}{3} & \color{g}{4} & \color{g}{3} & \color{g}{4} & \color{g}{4} & \color{g}{4} \\
 \hline
6&+&+&0.5&  & \bf1&  &  &  & \bf1&  & \bf1& \bf1& \bf1&  & \bf1& \bf1& \color{g}{2} & \bf1& \bf1& \bf1& \color{g}{2} & \bf1& \color{g}{2} & \color{g}{2} & \color{g}{2} & \bf1& \color{g}{2} & \color{g}{2} \\
6&+&-&0.5 &   &  &  &  & \bf1& \bf1&  &  & \bf1& \bf1& \bf1& \bf1& \bf1& \bf1& \bf1& \bf1& \color{g}{2} & \color{g}{2} & \bf1& \bf1& \color{g}{2} & \color{g}{2} & \color{g}{2} & \color{g}{2} & \color{g}{2} \\
6&-&+&0.5 &   &  & \bf1&  &  &  & \bf1& \bf1& \bf1&  & \bf1& \bf1& \bf1& \bf1& \color{g}{2} & \bf1& \bf1& \bf1& \color{g}{2} & \color{g}{2} & \color{g}{2} & \bf1& \color{g}{2} & \color{g}{2} & \color{g}{2} \\
6&-&-&0.5&   &  &  & \bf1&  & \bf1&  & \bf1&  & \bf1& \bf1& \color{g}{2} &  & \bf1& \bf1& \color{g}{2} & \bf1& \color{g}{2} & \bf1& \color{g}{2} & \bf1& \color{g}{2} & \color{g}{2} & \color{g}{3} & \bf1\\
 \hline
8&+&&1.5&   & \bf1&  & \bf1& \bf1& \color{g}{2} & \bf1& \color{g}{2} & \color{g}{2} & \color{g}{3} & \color{g}{2} & \color{g}{3} & \color{g}{3} & \color{g}{4} & \color{g}{3} & \color{g}{4} & \color{g}{4} & \color{g}{5} & \color{g}{4} & \color{g}{5} & \color{g}{5} & \color{g}{6} & \color{g}{5} & \color{g}{6} & \color{g}{6} \\
8&-&&1.5 &   &  & \bf1& \bf1& \bf1& \bf1& \color{g}{2} & \bf2 & \color{g}{2} & \color{g}{2} & \color{g}{3} & \color{g}{3} & \color{g}{3} & \color{g}{3} & \color{g}{4} & \color{g}{4} & \color{g}{4} & \color{g}{4} & \color{g}{5} & \color{g}{5} & \color{g}{5} & \color{g}{5} & \color{g}{6} & \color{g}{6} & \color{g}{6} \\
 \hline
\end{array}
}
\]
}
\caption{\label{dimensions} Dimensions of Atkin-Lehner subspaces of $S_k(N)^\epsilon$, with sources of rational newforms  highlighted in bold.}
\end{table}
the result was a factorizing quadratic.  On the table, $q$ is either $2$, $4$, 
or $8$.   

 The column $m$ gives the mass belonging to $(N,\epsilon)$,
calculated as a product of local masses \eqref{localmasses}.  The absence of a line
for $(N,\epsilon_4) =(4,+)$  reflects the vanishing $m^+(2,2)=0$ mentioned after 
\eqref{localmasses}.  The many ${\bf 1}$'s on the lines for $N=2$ arise from the 
small value $m^\pm(2,1)=1/2$.  The similarly many ${\bf 1}$'s for each line
belonging to $N=6$ arise because $m^{\pm}(2,1) = 1/2$ is not increased
by the second factor $m^{\pm}(3,1) = 1$.  

Table~\ref{dimensions}
clarifies the lines corresponding to $N=2$, $3$, $4$, $6$, and $8$ of Table~\ref{maintable}.
More importantly, it serves as an overview of the newforms studied in
 the next two sections.

\subsection{Sample computations at larger level} 
\label{SampleLarger}
To represent our computations for larger level,  
consider the weights $10$, $12$, $14$, and $16$ just 
below the weights $k \geq 18$ where Conjecture~\ref{mc} becomes effective. 
For these weights, the sequences of observed minimal levels, with
multiplicities, end as follows:    
\begin{equation}
\label{sequenceends} 
\begin{array}{ll}
k=10: & \dots,  210, {\bf 210}, {\bf 210}, {\bf 285}, {\bf 294}, {\bf 330}; \\
k=12: & \dots,  90, {\bf 96}, {\bf 114}; \\
k=14: & \dots, 42, {\bf 60}; \\
k=16: & \dots, 42, 42, {\bf 42}, {\bf 42}, {\bf 150}. \\
\end{array}
\end{equation}
We discussed the case $k=16$ in \S\ref{generalN} and \S\ref{dichotomy}.  The other
cases are similar, with ordinary type indicating forced
rationality and boldface unforced rationality.     For example, the forced $210$ comes from 
the one-element set $P_{10}(210)^{----}$, while the unforced $\bf{210}$, $\bf {210}$ comes
from the two-element set $P_{10}(210)^{+--+}$.    The largest set $P_k(N)^\delta$ 
where we have observed splitting in $k \geq 10$ is $P_{10}(285)^{-+-}$, which splits as 
$1+12$.

\begin{figure}[htb]
\begin{center}
\includegraphics[width=4.8in]{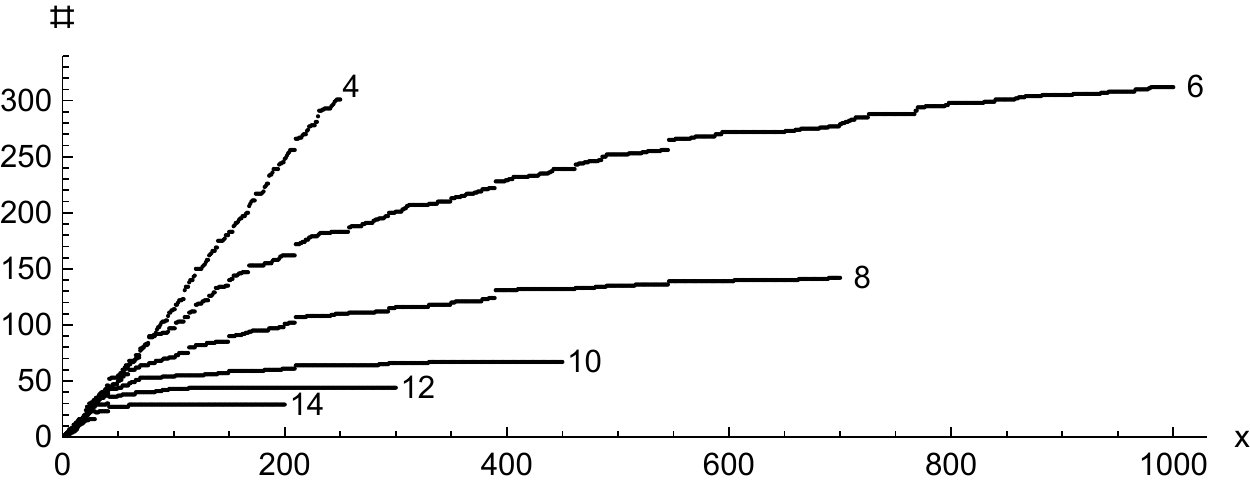}
\end{center}
\caption{\label{summatory} Graphs of the summatory functions $\#_k(x) = \sum_{N\leq x} |Q_k(N)|$,
for $k = 4$, $6$, $8$, $10$, $12$, and $14$.}
\end{figure}
Figure~\ref{summatory} graphs the summatory functions indicated by its caption.    
One has $\#_2(1000) = 1612$ and $\#_4(1000) = 802$, the cases $k=2$ and $k=4$
being not drawn and partially drawn respectively.  For $k = 6$ and $k=8$, the graphs
are still slightly rising as $x$ reaches the cutoff $C_k$. 
 For $k=10$, $12$, and $14$, the complete flattening
of the graph makes it plausible that the last minimal level seen $N_k'$ is indeed the last minimal 
level $N_k$.  The undrawn graph for $k=16$ almost coincides with the drawn
graph for $k=14$,  
although for $k=16$ the last minimal level seen $N_{16}'=150$ 
is close to the cutoff $C_{16}=200$.

 \section{Rings $M(N)$ for small $N$} 
 \label{rings}
 According to Conjecture~\ref{mc},  non-CM newforms in large weight 
 $k$ rarely have rational coefficients.   Our viewpoint is that those newforms which
 do have rational coefficients are of particular interest and deserve
 to be exhibited explicitly, in the style of Ramanujan's formulas \eqref{raman}.  
 We present some such explicit formulas here, systematically working
 in analogs of the classical ring $M(1) = \C[Q,R]$.  
 As a general convention, if $P_k(N)^\epsilon$ has just one
 element, then we call it $\Delta_{k,N}^\epsilon$.  
 
 We treat the cases $N = 2$, $3$, $4$, $6$, and $8$.  The new rings
 $M(N)$ contain $M(1)$ with indices $3$, $4$, $6$, $12$, and $12$.  
 We deal with the greater complexity by first embedding these rings
 into yet larger rings $\widehat{M}(N)$, all of which,
 like $M(1)$, are free on two generators.   We then keep formulas
 concise by systematically exploiting Atkin-Lehner operators in the 
 cases $N=2$, $3$, $6$, and $8$ and by discarding old forms cleanly
 in the cases $N=4$ and $8$.   
 
 An overall theme is that the five cases we treat are 
 remarkably similar to the classical case $M(1)$.    
 For example, the cuspidal ideal $S(N)$
 is always generated by an $\eta$-product
 analogous to the generator $\Delta_{12,1} = \eta^{24}$ of $S(1)$.   
 As analogs of Ramanujan's six newforms, Table~\ref{dimensions}
 for $N=2$, $3$, $4$, $6$, and $8$ lists 
 $23$, $13$, $6$, $47$, and $10$ newforms. 
 We keep things relatively brief by giving
 explicit formulas only for a small subset
 of these newforms.     Proofs of all statements
 are straightforward and omitted.  To
 proceed similarly for other $N$, 
 {\em Magma}'s \verb@Relations@ command 
 would be very useful.

 \subsection{Theta series}
 \label{level1}
 We build all our modular forms from two types of theta series,
 \begin{align*}
 \Theta & =   \sum_{(x,y) \in \Z^2} q^{x^2+x y+y^2} \!\!\!\!\!\!\!\!\!\!\!  \!\!\!\!\!\!\!\!\!\!\! && = 1 + 6 q + 6 q^3  + 6 q^4  + 12 q^7 + \cdots,  \\
 \theta & =    \sum_{x \in \Z} q^{x^2}  && = 1 + 2 q + 2 q^4 + 2 q^9 + 2 x^{16} + \cdots.
  \end{align*}
  These theta series are modular forms in their own right, with weights $1$ and $1/2$ respectively, and 
  certain characters.  Using standard notation, $\Theta \in M_1(3,\chi_{-3})$ and $\theta \in M_{1/2}(4,\chi_1)$. 
  As mentioned in \S\ref{genchar}, these two forms are outside of the context of the main review in Section~\ref{review}.  However, we
  are using them only to build forms which are in that context.  
  
  Using the second push-up operator $q \mapsto q^2$ of \S\ref{operators}, which does not change the weight, we obtain the graded rings 
  \begin{align*}
  \widehat{M}(6) & = \C[\Theta_1,\Theta_2],  & 
   \widehat{M}(8) & = \C[\theta_1,\theta_2].
   \end{align*}
   The rings of \S\ref{rings1} are then the subrings
    \begin{align*}
  M(6) & = \C[\Theta^2_1,\Theta_1 \Theta_2,\Theta^2_2],  & 
  M(8)  & = \C[\theta^4_1,\theta_1^2 \theta_2^2, \theta^4_2].
   \end{align*}
   These rings relate to their common subring $M(1)$ via
   \begin{align*}
   Q & =  \Theta_1 \left(5 \Theta_1^3+12 \Theta_1^2 \Theta_2-16 \Theta_2^3\right) \\ & 
    =  \theta_1^8+56 \theta_1^6 \theta_2^2-40 \theta_1^4 \theta_2^4-32 \theta_1^2 \theta_2^6+16 \theta_2^8  , \\
    &\\
    R & =  \left( \Theta_1^2+2 \Theta_1 \Theta_2-2 \Theta_2^2\right)
     \left( -11 \Theta_1^4-20 \Theta_1^3 \Theta_2+16 \Theta_1 \Theta_2^3+16 \Theta_2^4 \right)   \\
    & = \left(\theta_1^4+4 \theta_1^2 \theta_2^2-4 \theta_2^4\right) \left(\theta_1^8-136 \theta_1^6 \theta_2^2+152 \theta_1^4 \theta_2^4-32 \theta_1^2 \theta_2^6+16
    \theta_2^8\right).
   \end{align*}
   The four formulas just displayed are already a little bit long.  In our treatment of the various $M(N)$, we are
   exploiting structure to keep analogous formulas as short as possible.

   \subsection{Level $2$} 
   \label{level2} The ring $M(2)$ is freely generated by 
   \begin{align*}
   A & =  \theta_1^4+4 \theta_1^2 \theta_2^2-4 \theta_2^4  &\!\!\!\!\!\!\!\! & \in M_2(2)^- , \\
   B & = \theta_1^8-24 \theta_1^6 \theta_2^2+40 \theta_1^4 \theta_2^4-32 \theta_1^2 \theta_2^6+16 \theta_2^8 \!\!\!\!\!\!\!\! &\!\!\!\!\!\!\!\! & \in M_4(2)^-.
   \end{align*}
   The cuspidal ideal $S(2)$ is generated by 
   \[
   \Delta_{8,2}^+ = \eta_1^8 \eta_2^8 = 2^{-8} (A^4-B^2)   \in M_8(2)^+.
   \]
    Table~\ref{dimensions} shows $23$ rational newforms in $S(2)$, all of which have explicit expressions as polynomials
   in $A$ and $B$.  The one of largest weight is 
   \begin{eqnarray*}
   \Delta^{-}_{48,2} & = & 2^{-16} A^2 B  \left(49 A^4-81 B^2\right) \left(25 A^4-9 B^2\right) \cdot \\
   &&  \left(375531625 A^8-755257890 A^4 B^2+379726137
    B^4\right) \Delta^{+}_{8,2}.
    \end{eqnarray*}
    Here the fact that the total degree in $A$ and $B$ has to be odd cuts the number of terms roughly by half. 
    Throughout this section, there are many striking factorizations such as the one just displayed; 
    we are not pursuing their meaning.   
   
   \subsection{Level $3$} 
   \label{level3} 
   Here we work in the graded ring $\widehat{M}(3) = \C[\Theta,\Phi]$ with
   $\Phi = 4 \Theta_2^3 -3\Theta_1^2 \Theta_2 = 1-36q -54 q^2 -252 q^3 + \cdots \in M_3(3,\chi_{-3})$.  The even part of this
   ring is exactly $M(3)$.   The order four automorphism given by $\Theta \mapsto i \Theta$ and $\Phi \mapsto i \Phi$ restricts
   to the Atkin-Lehner involution $w_3$ on $M(3)$.   The cuspidal ideal $S(3)$ is generated by 
   \[
   \Delta_{6,3}^- =  \eta_1^6 \eta_3^6 = 2^{-2} 3^{-3} (\Theta^6-\Phi^2) = q - 6 q^2 + 9 q^3 + 4 q^4 + \cdots .
   \]
   The newform $\Delta_{6,3}^-$ is the first of the thirteen rational newforms with 
   level $3$ on Table~\ref{dimensions}.  The newforms from \S\ref{cardinality} exhibiting unforced 
   rationality are
   \begin{align*}
   \Delta_{22,3}^{+a} & =  6^{-1} \Theta \Phi  \; \;\; \; \left(75 \Theta^{12}-44 \Theta^6 \Phi^2-25 \Phi^4\right) \;\;\;\;\;\;  \Delta_{6,3}^- \;& \!\!\!\!&= q + 1728q^2 - 
      \cdots , \\
   \Delta_{22,3}^{+b} & =  3^{-3} \Theta \Phi  \left(-869 \Theta^{12}+1072 \Theta^6 \Phi^2-176 \Phi^4 \right) \Delta_{6,3}^- \;& \!\!\!\! & = q - 2844 q^2 - 
    \cdots .
   \end{align*}
   The difference of the coefficients of $q^2$ is $4572 = 2^2 3^2 127$.  Computation quickly suggests that
   in fact $\Delta_{22,3}^{+a} \equiv  \Delta_{22,3}^{+b}$.    In fact, seeing this congruence on the coefficient
   of $q^n$ for $n \leq 11$ suffices to confirm the general congruence, by a version of Sturm's theorem
   \cite[Prop.~1]{CKR10}.

   \subsection{Level $4$}  
   \label{level4}
    Here again we work in a larger graded ring $\widehat{M}(4) = \C[\theta,D]$, with 
   $\theta$ in weight $1/2$ and $D = \theta_1^4 - 8\theta_1^2\theta_2^2 + 8\theta_2^4 = 1-24 q + 24 q^2 - \cdots$ 
   in weight $2$.    The graded ring $M(4)$ is then just the sum of the graded pieces indexed by even integers.  
   Let $\rho$ be an eighth root of unity.  Then the automorphism $\theta \mapsto \rho \theta, D \mapsto D$ 
   restricts to the Atkin-Lehner involution $w_4$ on $M(4)$.    
   One thus has
   $M(4)=\C[C,D]$ with $C := \theta^4 \in M_2(4)^-$ and $D \in M_2(4)^+$.   
   Note that our presentation is a variant of \cite[IV.1 Prop. 4]{Kob93}, which uses the
   generators $\theta$ and $F = 2^{-5}(\theta^4-D)$.  
   
   For $k$ an odd integer, one has $\widehat{M}_k(4) = M_k(4,\chi_{-4})$.   For example, 
   the element $\Delta_{5,4} = 2^{-6} \theta^2 ( \theta^8+D^2) =  \eta_1^4 \eta_2^2 \eta_4^4 = q - 4 q^2 + 16 q^4 - \cdots$
   of $\widehat{M}_5(4)$ is a cuspidal newform with character $\chi_{-4}$ and CM.  
   Multiplying this element by $\theta^2$ we get a non-CM newform which fits into our 
   framework,
   \begin{equation}
   \Delta_{6,4}^- = \eta_2^{12} = 2^{-6} C (C^2-D^2) = q - 12 q^3 + 54 q^5 - 88 q^7 - \cdots .
   \end{equation}
   Thus one has the unusual situation of two quite different eta-products with quotient
   just $\theta^2$.  
   
   The $N=4$ block of Table~\ref{dimensions} has two interesting features.  First,
   as commented already  after \eqref{localmasses}, the $w_4=+$ part is zero and thus 
   missing from the block.  But second, the entries on the $w_4=-$ line 
   are exactly those of the familiar $N=1$ line, shifted to the left by a weight difference of six. 
    In fact, newforms on $N=4$ nicely separate
   from old forms via the formula
   \begin{equation}
   \label{nice4}
   S^{\rm new}_k(4)^{-} = M_{k-6}(1)_2 \Delta_{6,4}^- .
   \end{equation}
   The source of this equation is that the group $S_3$ acts on the 
   graded ring $M(4)$ with quotient $M(1)_2$, with the isotypical
   component corresponding to the sign character of $S_3$ being exactly $S^{\rm new}(4)^-$.  
   
   To be more explicit about \eqref{nice4}, one has $M(1)_2 = \C[Q_2,R_2]$ with $Q_2 = 2^{-2} (3 \theta^8 + D^2)$ and
   $R_2 = 2^{-3} D (9 \theta^8 - D^2)$.  Ramanujan's six rational forms \eqref{raman} in $S(1)$
   becomes six rational forms in $S^{\rm new}(4)^-$ via the simultaneous replacements $\Delta_{12,1} \mapsto \Delta_{6,4}^-$, 
   $Q \mapsto Q_2$, and $R \mapsto R_2$.   For example, the largest weight rational form on Table~\ref{dimensions} is
   $\Delta_{20,4}^- = Q_2^2 R_2 \Delta_{6,4}^-$.   
   
   \subsection{Level $6$}  
   \label{level6}
   The ring $\widehat{M}(6) = \C[\Theta_1,\Theta_2]$ and its even weight subring $M(6)$ have
   already been introduced in \S\ref{level1}.   The space $M_2(6)$ has dimension three, 
   with a basis consisting of   
   a sum, product, and difference: 
   \begin{align*} 
   s & \! = \! \Theta_1^2 \!+ \! 2 \Theta^2_2 \! \in \! M_2(6)^{+-}, &
   p & \! = \! \Theta_1 \Theta_2 \!  \in \! M_2(6)^{-+}, & 
   d & \! = \! \Theta_1^2 \! - \! 2 \Theta^2_2 \! \in \! M_2(6)^{--}.
   \end{align*}
   Via the equation $d^2-s^2=8p^2$, every element in $M(6)$ can be written in the canonical form $f_1(s,p)+d f_2(s,p)$.

      The cuspidal ideal $S(6)$ is generated by 
   \[
   \Delta_{4,6}^{++} = \eta_1^2 \eta_2^2 \eta_3^2 \eta_6^2 = 2^{-2} 3^{-2} (9 p^2-s^2).  
   \]
   Illustrations of how $\Delta_{4,6}^{++}$ is indeed a generator include
   \begin{align*}
   \Delta_{6,6}^{-+} & = p \Delta_{4,6}^{++}, &
   \Delta_{8,6}^{--} & = s p \Delta_{4,6}^{++}, &
   \Delta_{10,6}^{+-} & = 2^{-1} 3^{-2} (5 s^2 - 39 p^2) s \Delta_{4,6}^{++}.
   \end{align*}
   None of these expressions involve the generator $d$.  However the canonical expression
   for the rational newform of largest weight does:
   \begin{eqnarray*}
   \lefteqn{\Delta_{50,6}^{--} = 2^{-2} 3^{-9} d \cdot } \\
   && \left( 140349306081007255050000 p^{22}-111659120501660492670000 p^{20} s^2 
   \right. \\ && \left. 
   +27589681151783316300150 p^{18} s^4+2577120214736187574830 p^{16} s^6
   \right. \\ && \left. 
   -3234565067472714047760 p^{14} s^8+921682623552505460496 p^{12} s^{10} 
   \right. \\ && \left. 
   -149165289449290130931 p^{10} s^{12}+15554206382841117045 p^8 s^{14} 
   \right. \\ && \left. 
   -1070217851875219680 p^6 s^{16}+47245789680492400 p^4  s^{18}
   \right. \\ && \left. 
   -1218365734678125 p^2 s^{20}+14004203846875 s^{22} \right) \Delta_{4,6}^{++}.
  \end{eqnarray*}
   In comparison with the newform of the next largest weight, 
    $\Delta^-_{48,2}$ from \S\ref{level2},
   the expression for $\Delta_{50,6}^{--}$ is much longer.  Like for the case $N=2$, 
   our treatment of $N=6$ fully exploits the Atkin-Lehner operators, but does not cleanly
   discard old forms.   The difference in complexity can be attributed to the fact that 
   asymptotically $1/3$ of $\dim(S_k(2)^-)$ comes from newforms but only
   $1/6$ of $\dim(S_k(6)^{--})$ does.   The complexity of the displayed expression
   underscores the usefulness of discarding old forms cleanly, as 
   in \eqref{nice4} for $N=4$ and \eqref{nice8} for $N=8$.

   \subsection{Level $8$}  The ring $\widehat{M}(8) = \C[\theta_1,\theta_2]$ has already been described 
   in \S\ref{level1}.   Its integral weight part is generated by three forms in weight one, 
   $\theta_1^2, \theta_2^2 \in M_1(8,\chi_{-4})$ and $\theta_1 \theta_2 \in M_1(8,\chi_{-8})$.  
   The ring $M(8)=\sum_k M_k(8,\chi_{1})$ constitutes half of 
   the even weight part of $\widehat{M}(8)$.  A monomial $\theta_1^i \theta_2^j$ with
   total weight $(i+j)/2$ is in $M(8)$ if both $i$ and $j$ are even,
   and in the other half $\sum_k M_k(8,\chi_{8})$ if both $i$ and $j$ are odd.  
  
  The cuspidal ideal $S(8)$ of $M(8)$ is generated by 
   \[
   \Delta_{8,4}^+ = \eta_2^4 \eta_4^4 =  2^{-2} \theta_1^2 \theta_2^2 (-\theta_1^2+\theta_2^2) (\theta_1^2 - 2\theta_2^2) = q - 4 q^3 - 2 q^5 + 24q^7 - \cdots.
   \]
   Asymptotically, $S_k^{\rm new}(8)$ has one-quarter the dimension of $S_k(8)$.  Analogously to \eqref{nice4},
   each Atkin-Lehner eigenspace on newforms can be conveniently isolated via 
   \begin{equation}
   \label{nice8}
   S_k^{\rm new}(8)^\epsilon =  M(2)_2^\epsilon  \Delta_{8,4}^+.
   \end{equation}
   So all the ${\bf 1}$'s in the $N=8$ block of Table~\ref{dimensions} have monomial formulas of the form $E \Delta_{8,4}^{+}$, 
   with $E$ as follows:
   \[
   \begin{array}{c|cccccc}
   & 4 & 6 & 8 & 10 & 12 & 14 \\
   \hline
  Q_k(8)^+ & 1 && A_2^2 & A_2 B_2 & & A_2^3 B_2 \\
  Q_k(8)^-   & & A_2 & B_2 & A_2^3 & A_2^2 B_2 & 
   \end{array}.
   \]
   The unforced splitting in weight $k=16$ yielding a ${\bf \underline{2}}$ on Table~\ref{maintable} and a ${\bf 2}$ on Table~\ref{dimensions} is given by 
   \begin{align*}
   \Delta_{16,8}^{-a} & = 2^{-1} B_2 (23 A_2^4 - 21 B_2^2) \Delta_{8,4}^+  \!\!\!\!\!\!&& = q + 2700 q^3 - 251890 q^5 + \cdots, \\ 
   \Delta_{16,8}^{-b} & = 2^{-1} B_2 (-25 A_2^4 + 27 B_2^2) \Delta_{8,4}^+ \!\!\!\!\!\!&& = q - 3444 q^3 + 313358 q^5 -\cdots.
   \end{align*}
   Applying \cite[Prop.~1]{CKR10} as we did at the end of \S\ref{level3}, 
   one has $\Delta_{16,8}^{-a} \equiv \Delta_{16,8}^{-b}$ modulo $6144=2^{11} 3$.  
   
\section{Reductions modulo $\ell$ and associated number fields} 
\label{reps}
In general, the arithmetic of  coefficients
of  newforms is governed by Galois representations which in turn can be described
in terms of  number fields.  In particular, for 
a rational newform $\sum a_n q^n$, the reductions $a_n \in \Z/\ell^e$  are determined by a number field
with Galois group inside $GL_2(\Z/\ell^e)$.   In this section, we briefly discuss our 
ninety-nine newforms from this point of view.    Our goal is to make clear that 
the newforms here are a rich source of examples, 
in a way which complements the much studied $N=1$ 
case \cite{SwD73,Bos11}.

\subsection{Overview} 
 In the interest of 
brevity, we restrict to the cases $\ell^e \in \{2,3,5,7\}$.   
\begin{table}[htb]
\[
\begin{array}{c|cccccc|cccccc}
&\multicolumn{6}{c|}{\mbox{number of fields}} & 
\multicolumn{6}{c}{\mbox{number of newforms governed}} \\
\ell \setminus N & 1 & 2 & 3 & 4 & 6 & 8  & 1 & 2 & 3 & 4 & 6 & 8 \\
\hline
2 &  &&&&&&&&&&\\
3 &  &&&&&1&&&&&&10\\
5 &  &&1&1&1&2 &     &   & 7 & 6 & 24 & 7,3 \\ 
7 &  &1&2&1&2&4 &     &  12 & 8,3 & 6 & 24,12 & 2,4,2,2 \\
\hline
\multicolumn{7}{r}{\mbox{Total number of forms:}} &6 & 23 & 13 & 6 & 47 & 10 
\end{array}
\]
\caption{\label{pgl2} Summary of examples of number fields governing newforms}
\end{table}
For further brevity, we focus on number fields with Galois
group surjecting onto $PGL_2(\F_\ell)$.  We call the
other cases degenerate.    The numbers of 
$PGL_2(\F_\ell)$ number fields involved in various parts of this
section, blanks indicating $0$, are indicated in the
left half of Table~\ref{pgl2}.
Thus the classical case $N=1$ gives no examples
in our restricted context.   This degeneracy continues
somewhat into our setting $N \in \{2,3,4,6,8\}$.  
However as $\ell N$ increases, number fields with Galois group
all of $PGL_2(\F_\ell)$ begin to appear.
In \S\ref{defpoly} we give defining polynomials in all 
cases.   

A remarkable feature of our collection
of examples, indicated in the right half of Table~\ref{pgl2},
is that each of the sixteen number fields governs more 
than one rational newform.  The numbers listed
in the right half follow the alphabetic order of
the labeling in \S\ref{defpoly}.  For example, of the thirteen newforms
for $N=3$,  eight and three are governed modulo $7$ by $F_{3a}$ and
$F_{3b}$ respectively; the remaining two 
newforms are thus degenerate modulo $7$.  

\subsection{The projective correspondence for $\ell \leq 7$}
\label{matching}
Knowledge of the
$a_p$ for a newform  $\sum a_n q^n \in S^{\rm new}_k(N)$, for $p$ running over prime numbers, determines the entire 
newform via
the explicit formula
\begin{equation}
\label{L}
\sum_{n=1} \frac{a_n}{n^s} = \prod_{p|N} \frac{1}{1 - a_p p^{-s}} \cdot \prod_{p \nmid N} \frac{1}{1 - a_p p^{-s} + p^{k-1-s}}.
\end{equation}
Accordingly, we focus attention on the $a_p$.   When studying 
the reductions of the $a_p$ to $\F_\ell$, we exclude the
$p$ dividing $N \ell$, which 
behave differently.

The $a_p$ are governed by a Galois representation into $GL_2(\F_\ell)$ which 
we take to be semisimple, making it well-defined.  
To simplify, we will be working mainly not with  the $a_p$ themselves, rather with their normalized squares 
$s_p = a_p^2/p^{k-1} \in \Q$.    Let $f(x) \in \Z[x]$ be a degree
$\ell+1$ polynomial capturing the associated Galois representation
into $PGL_2(\F_\ell)$.  
Then the $s_p$, considered in $\F_\ell$, 
 are correlated with 
the partitions $\lambda_p$ giving the degrees of the irreducible
factors of $f(x)$ over $\Q_p$.  For $\ell \leq 7$, these correlations are as follows:
\begin{equation*}
{\renewcommand{\arraycolsep}{4pt}
\begin{array}{c|cc|ccc|ccccc|ccccccc}
& \multicolumn{2}{c|}{\ell=2} &  \multicolumn{3}{c|}{\ell=3} &  \multicolumn{5}{c|}{\ell=5} &  \multicolumn{7}{c}{\ell=7}  \\
\hline
s_p & 0 & 1 & 0 & 1 & 2 & 0 & 1 & 2 & 3 & 4 &  0 & 1 & 2 & 3 & 4 & 5 & 6 \\
\hline
\lambda_p & 21&3&  22 &31&4&    222&33&411&6&51&   2222 & 3311 & 44 &  611  & 71 & 8 & 8 \\ 
& 1^3&&  211 &1^4&4&    2211&&&6&1^6&   22211 & &  &   & 1^8 &  \\ 
\end{array}.
}
\end{equation*} 
The Galois group of $f(x)$ is all of $PGL_2(\F_\ell)$ if and only if 
all $s_p$ arise.  

In our restricted setting of ninety-nine newforms and four residual primes, a mod $\ell$ projective representation 
is either  surjective or has cyclic image.    In the latter case, the newform moreover satisfies a general
congruence of the form 
\begin{equation}
\label{eisenstein}
a_p \equiv p^i + p^j \; (\ell).
\end{equation}
Thus our restriction
to surjective representations corresponds to concentrating on the more mysterious
cases.  

\subsection{Explicit polynomials}
\label{defpoly}
For any of the ninety-nine newforms, all of the
$a_p$ with $p \nmid 2N$ are even.  Thus $s_p$ is always $0$ in $\F_2$ and the projective
Galois representations into $PGL_2(\F_2) = S_3$ are all
non-surjective.  For $\ell=3$, all 
newforms with $N < 8$ are likewise degenerate, but the ten newforms at $N=8$ 
are all nondegenerate and governed by 
\begin{align*}
\phi_{8}(x) & = x^4 - 2x^3 - 6x + 3, & \delta_8 & =     - 2^4 3^5.
\end{align*}
Here and later, next to each displayed polynomial $f(x)$ we show also the 
discriminant of the number field $\Q[x]/f(x)$.  In fact, for a given $k$, $N$, and 
$\ell$, theory, as partially summarized in \S\ref{discs} below, 
 gives only a small list of possibilities for these field discriminants.
The polynomials we display were found on the 
database \cite{jr-global-database}.  Matching as in \S\ref{matching} for all small
$p$ makes it very likely that the $f(x)$ are correct.  
For $\ell \leq 5$, the completeness of the database 
confirms correctness.  For all $\ell$, correctness
is confirmed using the Serre conjecture \cite{se,kw1,kw2} which implies that any 
$PGL_2(\F_\ell)$ field which is not totally real will
appear already in weight $k \leq \ell+1$.

For $\ell = 5$, Table~\ref{pgl2} says there are five fields.  Indexing by the relevant level as we did before, 
defining polynomials with Galois group $PGL_2(\F_5)$ and 
the indicated discriminant are 
\begin{align*}
f_{3}(x) & = x^6 - x^5 + 5x^4 - 5x^2 + 16x - 1, & d_{3} & = 3^4 5^9,  \\ 
f_{4}(x) & = x^6 - x^5 + 5x^3 + 10x^2 - 27x - 23, & d_4 &= 2^4 5^9, \\
f_{6}(x) & = x^6 - x^5 + 30x^3 - 15x^2 + 3x + 222, & d_{6} & = 2^4 3^4 5^9, \\
f_{8a}(x) & = x^6 - 2x^5 - 8x - 4, &  d_{8a} & = 2^6 5^7, \\
f_{8b}(x) & = x^6 - 2x^5 + 10x + 5,&  d_{8b} & = 2^6 5^9.
\end{align*}
One can sometimes describe the situation much more completely while 
still remaining brief.   For example at $N=6$,  the polynomial $f_6$
governs all twenty-four forms with $\epsilon_2 \epsilon_3 = -1$, while
all twenty-three forms with $\epsilon_2 \epsilon_3=1$ are degenerate. 

For $\ell = 7$, Table~\ref{pgl2} says there are ten fields.  Defining polynomials with Galois group $PGL_2(\F_7)$ and 
the indicated discriminant are 
\begin{align*}
F_{2}(x) & = x^8 - x^7 - 196 x^2 + 28 x - 28, &  D_{2} & = -2^6 7^{13}, \\
F_{3a}(x) & = x^8 - 4x^7 + 21x^4 - 21x^2 - 15x - 3,  & D_{3a} & = -3^6 7^{11}, \\
F_{3b}(x) & = x^8 - 3x^7 - 7x^6 + 49x^5 + 42x^4,&  D_{3b} & = - 3^6 7^{13}, \\
F_{4}(x) & = x^8 - x^7 - 7x^6 + 7x^2 - 27x - 1, & D_4 & = -2^4 7^{11}, \\
F_{6a}(x) & = x^8 - 2x^7 + 42x^4 - 126x^3 + 84x^2 + 66x - 48, & D_{6a} & = -2^6 3^6 7^9, \\
F_{6b}(x) & =  x^8 - x^7 + 21x^6 + 21x^5 - 21x^4 + 945x^3 & D_{7b} & = -2^6 3^6 7^{13}, \\ 
& \qquad - 441x^2 + 45x + 3168, \\
F_{8a}(x) & =  x^8-2 x^7+7 x^4-14 x^2+8 x+5, & D_{8a}& = - 2^8 7^9,   \\
F_{8b}(x) & =  x^8-2 x^7-14 x^4+28 x^2-60 x+92, & D_{8b}&  = - 2^8 7^{11}, \\
F_{8c}(x) & = x^8-2 x^7+14 x^6+42 x^5+140 x^4 &  D_{8c}& = - 2^8 7^{13},   \\
& \qquad  +266 x^3+322 x^2+222 x - 157, \\
F_{8d}(x) & =  x^8-2 x^7+49 x^4-196 x^2-140 x-63, & D_{8d}& = - 2^8 7^{13}.  
\end{align*}
All twelve newforms in $M(2)^+$
are governed by $F_2$, while the eleven newforms in $M(2)^-$ are all degenerate.
The twenty-four newforms in $M(6)$ with $\epsilon_2=+$ are governed by
$F_{6a}$; of those in $M(6)$ with $\epsilon_2=-$, twelve are governed by $F_{6b}$ and
eleven are degenerate.

\subsection{Field discriminants and ramification}
\label{discs} General facts, some used to find the polynomials of the previous
subsection, are  illustrated by the displayed discriminants $D$.   Certainly, all primes dividing the discriminant
must divide $N \ell$ and for odd $\ell$ the discriminant must be a square times $\chi_{-4}(\ell) \ell$.  

Ramification at $\ell$ is directly related to weight.  
As mentioned earlier, all  polynomials necessarily arise already from newforms in weight $k \leq \ell+1$.  
In fact if $\ord_\ell(D) \geq \ell+2$, then in this range the polynomial arises only 
in weight $\ord_\ell(D) + 2- \ell$.   Thus in the reduction bijection from 
$\{\Delta_{4,8}^+,\Delta_{6,8}^-, \Delta_{8,8}^+, \Delta_{8,8}^-\}$
to $\{ f_{8a},f_{8b},f_{8c}, f_{8d}\}$, discriminants force the correspondences  
$\Delta_{4,8}^+ \leftrightarrow f_{8a}$ and $\Delta_{6,8}^- \leftrightarrow f_{8b}$. 
In fact, the only ambiguity as to the
canonical lowest weight source of each of our fifteen polynomials 
is the rest of this bijection.  For $\Delta_{8,8}^+$ and $\Delta_{8,8}^-$, one has 
$s_3 = 6$ and $0$ respectively.  For $f_{8c}$ and $f_{8d}$, one has $\lambda_3 = 8$ and $22211$ 
respectively.  By \S\ref{matching}, the bijection is completed by 
$\Delta_{8,8}^+ \leftrightarrow f_{8c}$  and $\Delta_{8,8}^- \leftrightarrow f_{8d}$.

In general, ramification at primes $p$ different from $\ell$ is directly related to the refinement of the Atkin-Lehner decomposition
discussed in \S\ref{generalN}.   With our tiny levels, we see only a small part of this complicated theory.
If $p$ exactly divides $N$, then
$\ord_p(D)$ is usually $\ell-1$ but exceptionally can be zero.  For the ten instances in 
\S\ref{defpoly}, it is always $\ell-1$.    If $\ord_2(N)$ is $2$ or $3$, then the 
$2$-decomposition group has to be $S_3$ and $S_4$ respectively; these nonabelian subgroups
exclude the simple abelian behavior \eqref{eisenstein};  they 
partially explain why all newforms at these levels $4$ and $8$ have
surjective mod $5$ and mod $7$ projective representations.   For 
$\ord_2(N)=3$, the slope content as in \cite{jr-local-database} has
to be $[4/3,4/3]_3^2$, as opposed to the other possibility for $S_4$ $2$-adic fields,
$[8/3,8/3]_3^2$.  This restriction accounts for the small exponents on $2$ in the
seven polynomials in \S\ref{defpoly} associated to $N=8$.  

\subsection{Congruences}  A necessary condition for two rational newforms with the same level $N$
to reduce to the same power series in $\F_\ell[[q]]$ is that their projective mod $\ell$ representations 
coincide and their weights are congruent modulo $\ell-1$.   This condition is sufficient for $\ell=2$, but not 
for $\ell >2$, as there is still a sign ambiguity in each of the $a_p$.   For 
example, the two newforms $\Delta^{\pm}_{8,8}$ satisfy the necessary condition,
but differ via twisting by $\chi_{-3}$, as illustrated by 
$(a_5,a_7,a_{11},a_{13},a_{17},a_{19},a_{23},a_{29}) = (\pm 1, 0, \pm 1, 1, \pm 2, 2,\pm 1,\pm 0)$.

Remarkably, the sign ambiguity can be resolved in a simple way in all our cases.  Most strikingly, for $N \in \{2,4,6\}$
the above necessary condition is also sufficient, and in the degenerate case only the congruence
condition on weights needs to be verified.     As examples with $\ell=7$,
\begin{eqnarray*}
\Delta^+_{8,2} \equiv \Delta^+_{14,2} \equiv \Delta^+_{20,2} \equiv \Delta^+_{26,2} & \equiv  &
 q + 6q^2 + 5q^3 + q^4 + 2q^6 + q^7 + 6q^8 + \cdots, \\
\Delta^-_{10,2} \equiv \Delta^-_{22,2} \equiv \Delta^-_{28,2} \equiv \Delta^-_{40,2}  & \equiv&
 q + 2q^2 + 5q^3 + 4q^4 + 2q^5 + 3q^6 + q^8 + \cdots . 
\end{eqnarray*}
Other similar examples can be read off from Table~\ref{dimensions}.

\section{Concluding discussion}
\label{discussion}
     Sections~\ref{rings} and \ref{reps} examined particular non-CM 
     newforms with rational coefficients. 
      Here we return to 
   discussing  Conjecture~\ref{mc}, which concerns the landscape of all non-CM 
   newforms with rational
   coefficients.  Given that the Maeda conjecture has been open for twenty years, 
   we do not expect that
 our similar Conjecture~\ref{mc} will be proved soon.  In contrast, our assessment
 is that there is still insight to be gained by 
 continuing exploratory
 computations in this spirit of this paper.  Our concluding discussion
 proposes three directions for such computations.
 
 \subsection{Larger cutoffs.}  We have simply used {\em Magma}'s general modular form package to compute the sets $Q_k(N)$.  
 Programs optimized for this exact problem could likely go further, meaning larger cutoffs $C_k$ for each given $k$.   Steps have been 
 taken in this direction \cite{KedMed}, with one of the main ideas being to first work with modular forms modulo two.

 \subsection{Connections with extra vanishing}  
 \label{extravan} Let $P_k(N)^ \delta$ be a locally defined collection of newforms, 
 as in \S\ref{generalN}.  Every newform $g$ in this set has a completed $L$-function $\Lambda(g,s)$, with functional equation 
 $\Lambda(g,s) = w \Lambda(g,k-s)$.  The sign $w$ of the functional equation is the same for all 
 $g$.   As a consequence,
 the order of vanishing $r(g)$ at the central point $s = k/2$ satisfies $(-1)^{r(g)} = w$ and so has constant parity.  
 While the rank $r(g)$ itself can certainly vary with $g$,
 the Beilinson-Bloch arithmetic interpretation of central vanishing implies that conjugate $g$ should have identical $r(g)$.  
 
 Vanishing beyond order one therefore has the potential to ``explain'' some locally unforced splittings.  For example, we have looked at all 
 rational non-CM newforms $g$ with $k=6$ and $N \leq 400$.  
 Eight of them have $r(g) = 2$ and none
 have $r(g) \geq 3$.   Seven of these eight have quadfree level and are members of sets $P_6(N)^\epsilon$
 as follows:
\[
\begin{array}{c|ccccccc}
\multirow{2}{*}{$N$} & 95 & 116 &  122 &  260 &  308 &  359 & 371 \\
    & 5 \cdot 19 & 4 \cdot 29 & 2 \cdot 61 & 4 \cdot 5 \cdot 13 & 4 \cdot 7 \cdot 11 & 359 & 7 \cdot 53   \\
\hline
\epsilon & +- & -+ & +- & -++ & -++ & - & +- \\
d &         9+1 & 6+1 & 6+1 & 4+1 & 6+1 & 83+1 & 34+1  
\end{array}.
\]
In each case, we work with the smallest $p$ not dividing $N$.  The degree $d$ 
characteristic polynomial $f_{6,N,p}^\epsilon(x)$ always factors as a linear factor 
times a polynomial with Galois group $S_{d-1}$.   For example $f_{6,359,2}^-(x) = (x-5)(x^{83}- 7 x^{82} - \cdots)$ 
with the degree $83$ polynomial having a $4128$-digit field discriminant.

The order of vanishing $r(g^\chi)$ of quadratic twists $g^\chi$ can also be 
taken into consideration.   As an example, let $g= q + 4q^2 + 11 q^3 + \cdots$ be the unique form in the 
space denoted $P_6(50)^{-\pm}$ in \S\ref{generalN}.  Both $g$ and its twist $g^{\chi_5} = q - 4 q^2 - 11 q^3 + \cdots   \in P_6(50)^{+\mp}$ 
have rank zero.   However the twist $g^{\chi_{-4}}$ with level $400$ has rank two.   To deal 
with arbitrary quadratic twists of a fixed form $g$, it would be natural to bring in half-integral
weight modular forms via the Shimura-Waldspurger correspondence as described with
examples in \cite[IV.4]{Kob93}.    

In weight six, we computed also in levels larger than $400$ 
and have seen several more examples of 
extra vanishing, all again with rank two.    We have not seen any extra vanishing at all in weights $\geq 8$,
and so in particular we have no explanation of the various unforced splittings seen there.  However, 
we have not computed systematically and there may be extra vanishing for twists at larger levels.   
In weight eight, order two vanishing has been seen in CM newforms \cite[\S6.6.2]{Wat08}.

\subsection{A broader notion of rationality}  This paper has addressed the problem of tabulating all 
non-CM newforms which
satisfy the rationality condition that all their Fourier coefficients $a_n$ are rational.    There is second much
weaker rationality condition that is equally natural from a motivic point of view, 
namely that all the $|a_n|^2$ are rational.    In this larger setting, one has to work with general characters $\chi$ 
and allow odd weights $k$ as well.   The rings described in Section~\ref{rings}
  include some interesting
examples.  

The two problems are both instances of a common general problem
concerning objects in the category $\mathcal{M}(\Q,\Q)$ of 
motives over $\Q$ with coefficients in $\Q$.   Roughly
speaking, the problem of this paper is equivalent to 
classifying rank two motives in 
this category with Sato-Tate group 
the symplectic group $Sp_2$.  The larger 
problem with the weaker rationality condition is equivalent to 
classifying rank three motives with Sato-Tate group
the special orthogonal group $SO_3$.  With this weaker
notion of rationality, there 
would of course be more newforms to collect; however
we would still expect finiteness in all weights 
$k \geq 5$.

\bibliographystyle{alpha}
\bibliography{nfr}

 \end{document}